\documentclass[10pt]{article}

\usepackage{a4wide}
\usepackage{amssymb}
\usepackage{amsfonts}
\usepackage{amsmath}
\input xy
\xyoption{arrow} \xyoption{matrix}

\date{}

\newtheorem{proposition}{Proposition}[section]
\newtheorem{theorem}[proposition]{Theorem}
\newtheorem{lemma}[proposition]{Lemma}

\newtheorem{corollary}[proposition]{Corollary}

\def\der{\partial }

\def\nFM0{{\nu }_{F,M_0}}
\def\nFN0{{\nu }_{F,N_0}}
\def\nGN0{{\nu }_{G,N_0}}

\def\N0{ {\bf N}_0 }

\def\g{\gamma}

\def\ra{\rightarrow}

\def\lra{\leftrightarrow}
\def\Xpm{X^{\pm }}

\def\s{\sigma}
\def\Z{{\bf Z }}

\def\l1{{\lambda}_1}

\def\a{\alpha}
\def\a0{ {\alpha }_0}
\def\a1{ {\alpha }_1}

\def\l{\lambda}
\def\o{\omega}
\def\nFGM0{{\nu }_{F,G,M_0}}

\def\nFN0{{\nu}_{F,N_0}}

\def\sm{{\sigma}^m}

\def\sm1{{\sigma}^{-1}}

\def\smtp1{{\sigma}^{-t+1}}

\def\o{\omega }
\def\S1{S^{-1}}

\def\Xpm1{X^{\pm 1}_1}

\def\sPM1{{\sigma }^{\pm 1}}
\def\sMP1{{\sigma }^{\mp 1 }}

\def\b{\beta}
\def\d{\delta}

\def\di{{\rm d.ind}}

\def\L{\Lambda}

\def\OO{{\cal O}}

\def\CD{{\cal D}}

\def\Ytm1{Y^{t-1}}
\def\Yim1{Y^{i-1}}

\def\CK{{\cal K}}

\def\ass{{\rm ass}}

\def\Aut{{\rm Aut}}

\def\ker{ {\rm ker } }

\def\D{ \Delta }
\def\SL2Z{ {\rm SL}_2({\bf Z}) }

\def\th{ \theta }

\def\Gp1{ G^{1 , 1 } }
\def\P11{ P^{-1 , 1 } }
\def\Pp1{ P^{1 , 1 } }

\def\Supp{{\rm Supp}}

\def\th{\theta}

\def\nCLsr{{}^\nu\kern-2pt {\cal L}^{\sigma , \rho  }}
\def\nP{{}^\nu \kern-2pt P}
\def\nL{{}^\nu\kern-2pt L}
\def\nLL{{}^\nu\kern-2pt \Lambda}
\def\nPsr{{}^\nu\kern-2pt P^{\sigma , \rho  }}
\def\nLsr{{}^\nu\kern-2pt L^{\sigma , \rho  }}
\def\nuCL{{}^\nu\kern-2pt  {\cal L}}
\def\nCLsr{{}^\nu\kern-2pt {\cal L}^{\sigma , \rho  }}
\def\nCL1m{{}^\nu\kern-2pt {\cal L}^{-1 , 1  }}
\def\x1nu{x^\frac{1}{\nu}}
\def\xm1nu{x^{-\frac{1}{\nu}}}

\def\ra{\rightarrow }

\def\CB{{\cal B}}

\def\CC{ {\cal C}}

\def\CP{ {\cal P}}

\def\nAM0{{\nu }_{{\cal A},M_0}}
\def\nAN0{{\nu }_{{\cal A},N_0}}

\def\End{ {\rm End }}

\def\CP{ {\cal P }}

\def\ga{\mathfrak{a}}
\def\gb{\mathfrak{b}}

\def\gp{\mathfrak{p}}

\def\di!{\frac{\der^i}{i!}}
\def\dik!{\frac{\der^k_i}{k!}}

\def\CC{{\cal C}}

\def\id{{\rm id}}

\def\N{\mathbb{N}}

\def\0{\overline{0}}
\def\1{\overline{1}}

\def\Ln1{\L_{n,\overline{1}}}

\def\a1{a_{\overline{1}}}

\def\S{\Sigma}

\def\vn1{\overrightarrow{n-1}}

\def\Q{\mathbb{Q}}
\def\im{{\rm im}}

\def\sl{{\rm sl}}

\def\mL{\mathbb{L}}

\def\mJ{\mathbb{J}}
\def\mI{\mathbb{I}}

\def\mF{\mathbb{F}}

\def\K1{{\rm K}_1}

\def\Supp{{\rm Supp}}

\def\hmI1{\widehat{\mI_1}}
\def\tmI1{\widetilde{\mI_1}}
\def\tmJ1{\widetilde{\mJ_1}}
\def\hB1{\widehat{B_1}}
\def\hCB1{\widehat{\CB_1}}

\def\Den{{\rm Den}}

\def\Ore{{\rm Ore}}

\def\Den{{\rm Den}}

\def\ga{\mathfrak{a}}

\def\pCC{{}'\CC}

\def\Z{\mathbb{Z}}

\def\mL{\mathbb{L}}
\def\mR{\mathbb{R}}
\def\lCCD{{'\CC_D}}
\def\lCCA{{'\CC_A}}
\def\tn{\widetilde{n}}

\begin{document}

\author{V. V. \  Bavula %(GWA-di-skew.tex)
}

\title{Generalized Weyl algebras and diskew polynomial rings}

\maketitle
\begin{abstract}

The aim of the paper is to extend the class of generalized Weyl algebras to a larger class of rings (they are also called {\em generalized Weyl algebras})  that are determined by two ring  endomorphisms rather than one as in the case of `old' GWAs.  A new class of rings, the {\em diskew polynomial rings}, is introduced that is closely related to GWAs (they are GWAs under a mild condition). The, so-called, ambiskew polynomial rings are a small subclass of the class of diskew polynomial rings. Semisimplicity criteria are given for generalized Weyl algebras and diskew polynomial rings.

 {\em Mathematics subject classification 2010: 16D30, 16P40, 16D25, 16P50, 16S85.}

$${\bf Contents}$$
\begin{enumerate}
\item Introduction. \item  Generalized Weyl algebras. \item Simplicity criteria for Generalized Weyl algebras. \item  Generalized Weyl algebras of rank $n$.

\item Diskew polynomial rings.  
\end{enumerate}

\end{abstract}

%%%%%%%%%%%%%%%%%% SECTION 1 %%%%%%%%%%%%%%%%%%%%%%%%

\section{Introduction}

{\bf Generalized Weyl algebras $D(\s , a)$ with central element $a$}.

 {\it Definition}, \cite{Bav-GWA-FA-91}-\cite{Bav-Bav-Ideals-II-93}. Let $D$ be a ring, $\s $ be a ring
automorphism of $D$, $a$ is a {\em central} element of $D$. The
{\bf generalized Weyl algebra} of rank 1 (GWA, for short) $D(\s , a)=D[x,y; \s , a]$ is a
ring generated by the ring $D$  and two elements $x$ and $y$ that
are subject to the defining relations:
%\marginpar{clGWA}
\begin{equation}\label{clGWA}
 xd=\s (d) x\;\;  {\rm and} \;\;  yd=\s^{-1} (d)y\;\; {\rm for \; all} \;\; d\in D,  \;\;
 yx=a \;\; {\rm and} \;\;   xy=\s (a).
\end{equation}
The ring $D$ is called the {\em base ring} of the GWA. The automorphism $\s$ and the element $a$ are called the {\em
defining automorphism} and the {\em defining element} of the GWA, respectively.

 This is an experimental fact that many popular algebras of small
Gelfand-Kirillov dimension are GWAs (see Section \ref{GWA2}): the first Weyl algebra $A_1$
and its quantum analogue, the {\em quantum plane}, the {\em
quantum sphere}, $Usl(2)$, $U_qsl(2)$, the {\em Heizenberg}
algebra and its quantum analogues, the $2\times 2$ quantum
matrices, the {\em Witten's} and {\em Woronowic's} deformations,
Noetherian down-up algebras, etc.

The generalized Weyl algebras were introduced by myself in 1987  when I was an algebra postgradute student at the Kyiv University, the Department of Mathematics,  and they were the subject of my PhD ``Generalized Weyl algebras and their representations'' submitted at the end of 1990 (defended at the beginning of 1991).

The aim of the paper is to introduce a generalization of GWAs where the elements $x$ and $y$ act on the ring $D$ by two {\em non-commuting} (in general) ring endomorphisms $\s$ and $\tau$ and the element $a$ is not central but a left normal.

{\bf Generalized Weyl algebras with two endomorphisms and a left normal element $a$.}

{\it Definition}. Let $D$ be a ring, $\s $ and $\tau$ be ring endomorphisms of  $D$, and an element $a\in D$ be  such that
%\marginpar{abGWA}
\begin{equation}\label{abGWA}
\tau \s (a) = a, \;\; ad= \tau \s (d) a\;\; {\rm and}\;\; \s (a) d= \s\tau (d) \s (a) \;\; {\rm for \; all}\;\; d\in D.
\end{equation}
The {\bf generalized Weyl algebra} (GWA) of rank 1, $A= D(\s, \tau, a) = D[x,y; \s, \tau, a]$,  is a ring generated by $D$, $x$ and $y$ subject to the defining relations:
%\marginpar{GWADEF}
\begin{equation}\label{GWADEF}
xd=\s (d) x\;\;  {\rm and} \;\;  yd=\tau (d)y\;\; {\rm for \; all} \;\; d\in D,
 \;\; yx=a \;\; {\rm and} \;\;   xy=\s (a).
\end{equation}
The ring $D$ is called the {\em base ring} of the GWA $A$. The endomorphisms $\s$, $\tau$ and  the element $a$ are called the {\em
defining endomorphisms} and the  {\em defining element} of the GWA $A$, respectively. By (\ref{abGWA}), the elements $a$ and $\s (a)$ are left normal in $D$. An element $d$ of a ring $D$ is called {\em left normal} (resp., {\em normal}) if $dD\subseteq Dd$ (resp., $Dd=dD$). See Section \ref{GWA2} for more information about left/right normal elements. To distinguish `old' GWAs from the `new' ones the former are called the {\em classical} GWAs. Every classical GWA is a GWA as the conditions in (\ref{abGWA}) trivially hold if $a$ is central and $\tau = \s^{-1}$. 

The construction of GWAs with two endomorphisms existed in 1987 but there were no natural examples at the time which is not surprising in view of the conditions (\ref{abGWA}) the defining element $a$ must satisfy. `Artificial'  examples of GWAs with two endomorphisms can be easily constructed from a GWA $A= D[x,y; \s , a]$ with a single automorphism where $a\in Z(D)$ by replacing the elements $x$ and $y$ by $x'=ux$ and $y'=vy$ where $u$ and $v$ are noncentral units or left normal elements of $D$.

{\bf Simplicity criteria for generalized Weyl algebras.}
 Theorem \ref{1Apr} shows existence of generalized Weyl algebras.  Let $D$ be a ring and $\s$ be its ring endomorphism.  An ideal $I$ of $D$ is called $\s$-{\em stable} if $\s (I) \subseteq I$. The ring $D$ is called a $\s$-{\em simple} ring iff $0$ and $D$ are the only $\s$-stable ideals of the ring $D$. An endomorphism $\s$ is {\em inner} if $\s = \o_u$ for some unit $u\in D$ ($\s (d)= udu^{-1}$ for all $d\in D$). Then necessarily $\s$ an automorphism of $D$.  The results of  Section \ref{GWA2} are used in the proofs of simplicity criteria for generalized Weyl algebras in Section \ref{SIMCGWA} (Theorem \ref{B5Apr} and Theorem \ref{29Oct}). Theorem \ref{29Oct}
 is a simplicity criterion for an arbitrary GWA and Theorem \ref{B5Apr} is a simplicity criterion for a GWA with a mild restriction (the elements $a$ and $\s (a)$ are normal), its proof is given in Section \ref{SIMCGWA}. Theorem \ref{B5Apr} generalizes  a simplicity criterion for the  (classical) GWAs  {\cite[Theorem 4.2]{Bav-FilDimSimCrit-1996}}.
\begin{theorem}\label{B5Apr}%\marginpar{B5Apr}
 Let $A= D[x,y; \s, \tau , a]$ be a GWA such that  the elements $a$ and $\s (a)$ are right  normal in  $D$ (they are normal, by (\ref{abGWA})). Then the following statement are equivalent.
\begin{enumerate}
\item  $A$ is a simple ring.
\item
\begin{enumerate}
\item The elements $a$ and $\s (a)$ are regular in $D$,
\item $D$ is a $\s$-simple ring,
\item for all $i\geq 1$,   $\s^i$ is not an inner automorphism of the ring $D$, and
    \item for all $i\geq 1$, $Da+D\s^i(a) = D$.
\end{enumerate}
\item
\begin{enumerate}
\item The elements $a$ and $\s (a)$ are regular in $D$,
\item $D$ is a $\tau$-simple ring,
\item for all $i\geq 1$,   $\tau^i$ is not an inner automorphism of the ring $D$, and
    \item for all $i\geq 1$,  $D\s (a) + D\tau^i \s(a) = D$.
\end{enumerate}
\end{enumerate}
If one of the equivalent conditions holds then $\s$ and $\tau$ are automorphisms of $D$.
\end{theorem}
At the end of Section \ref{SIMCGWA}, natural classes of involutions on GWAs  are introduced (Lemma \ref{b2Apr}).

 {\bf Generalized Weyl algebras of arbitrary rank.}
In Section \ref{GWARNN}, a more general definition of generalized Weyl algebras of {\em arbitrary}  rank is given. The classical generalized Weyl algebras of rank $n$ are determined by $n$ commuting automorphisms and $n$ central elements of the base ring, \cite{Bav-GWA-FA-91}-\cite{Bav-Bav-Ideals-II-93}, see Section \ref{GWARNN}. The `new' GWAs are determined by $2n$ not necessarily commuting  ring endomorphisms and $n$ left normal  elements of the base ring. A large class of examples are considered (Proposition \ref{31Jul}).

{\bf Diskew polynomial rings}. Let us introduce a new class of rings -- the diskew polynomial rings (DPR).

{\it Definition}. Let $D$ be a ring, $\s$ and $\tau$  be its ring endomorphisms, $\rho$ and $b$  be elements of $D$ such that, for all $d\in D$,
%\marginpar{wbd2}
\begin{equation}\label{wbd2}
\s\tau (d)\rho = \rho \tau \s (d)\;\; {\rm and}\;\; \s \tau (d)b=bd,
\end{equation}
The {\bf diskew polynomial ring} (DPR)  $E:= D(\s , \tau, b , \rho ):= D[x,y; \s , \tau , b,\rho ]$ is a ring generated by $D$, $x$ and $y$ subject to the defining relations:
%\marginpar{DiskewDEF}
\begin{equation}\label{DiskewDEF}
xd=\s (d) x\;\; {\rm and}\;\; yd= \tau (d) y\;\; {\rm for\; all}\;\; d\in D, \;\; xy-\rho yx = b.
\end{equation}
By (\ref{wbd2}), $b$ is a left normal element of $D$. If $\tau\s$ (resp., $\s\tau$) is an epimorphism then $\rho$ is a left (resp., tight) normal element of $D$. Theorem \ref{5Apr} shows  existence of diskew polynomial rings.
The diskew polynomial rings are a generalization of the following class of rings which is a part of the class of diskew polynomial rings.

{\it Definition}, \cite{Bav-GlGWA-1996}. Let $D$ be an ring and $\s$ be its automorphism.  Suppose that elements $b$ and $\rho $
belong to the centre of the ring $D$,  $\rho $ is invertible and
 $\s (\rho )=\rho$.  Then   $E:=D\langle \s ;b,\rho\rangle:= D\langle X, Y; \s ,b,\rho \rangle$ is a ring generated by  $D$, $X$ and $Y$ subject to the defining
relations:
%\marginpar{Dsbrho}
\begin{equation}\label{Dsbrho}
X\alpha =\s (\alpha )X\;\; {\rm and} \;\; Y\alpha =\s {}^{-1}(\alpha )Y \;\; {\rm for \; all}\;\;   \alpha \in D, \;\; {\rm and} \;\;\, XY-\rho
YX=b.
\end{equation}
Clearly, $E=D[X,Y; \s  , \s^{-1} , b, \rho ]$ is a diskew polynomial ring. The origin of this construction stems from the universal enveloping algebra $Usl(2)$ of the Lie algebra $sl(2)$.  When we rewrite the defining relations of $Usl(2)$  (where $[a,b]=ab-ba$): $[H,X]=X$, $[H,Y]=-Y$ and $[X,Y]=2H$ in the equivalent form: $XH=(H-1)X$, $YH=(H+1)Y$ and $XY-YX=2H$ and notice that $XH= \s (H) X$ and $YH=\s^{-1}(H)Y$ where $\s$ is an automorphism of the polynomial algebra $K[H]$ given by $\s (H) = H-1$ we come to $Usl(2)=K[H]\langle X,Y; \s , 2H, 1\rangle$, \cite{Bav-GWA-Clas-sl2-90, Bav-GWA-FA-91}. The next natural step was to replace the polynomial $2H$ by an arbitrary polynomial $a(H)\in K[H]$. This was done independently in \cite{Bav-GWA-FA-91} and  \cite{Smith-SimUsl2-1990}. That is how the, so-called, {\em algebras similar to} $Usl(2)$  appeared. It is the algebra $K\langle X,Y,H\rangle$ that satisfies the defining relations:
$$XH=(H-1)X, \;\; YH=(H+1)Y\;\; {\rm  and}\;\; XY-YX=a(H) \;\; {\rm where}\;\; a(H)\in K[H]. $$
In 90s, there were many examples like this, various `quantum deformations' of $Usl(2)$, with a ring $D$ which is a `small' commutative ring, eg, $U_q(sl_2)$, $\OO_{q^2}(so(K,  3))$, the quantum Weyl algebra, the quantum plane, etc (see Section \ref{GWA2}).

If $D$ is commutative domain, $\rho =1$ and $b = u-\s (u)$ for some $u\in D$ (resp., if $D$ is a commutative finitely generated domain over a field $K$ and  $\rho \in K^*$) the algebras $E$ were considered in \cite{Jordan-ItSkewPol-1993} (resp., \cite{Jordan-Wells-1996}).

The ring $E=D\langle \s ;b,\rho\rangle$ is {\em the iterated skew polynomial ring}
$E=D[Y;\sm1  ][X;\s ,\der ]$ where $\der $ is the $\s -$derivation of
$D[Y;\sm1 ]$ such that $\der (D)=0$ and $\der (Y)=b$ (here the automorphism $\s $ is extended
from $D$ to $D[Y;\sm1 ]$ by the rule $\s (Y)=\rho Y)$.

{\bf Diskew polynomial rings are GWAs when $\rho$ is a unit}. If the element $\rho$ is a unit in $D$ then every diskew polynomial ring is a generalized Weyl algebra, Theorem \ref{2Apr} (a proof is given in Section \ref{DPR}).
\begin{theorem}\label{2Apr}%\marginpar{2Apr}
Let $E=D[x,y; \s , \tau , b, \rho ]$ be a diskew polynomial ring. Suppose that $\rho$ is a unit in $D$. Then   $x$ and $y$  are  left regular elements of $E$ and  the ring $E= \CD [ x,y; \s , \tau , a=h]$ is a GWA with base ring $\CD := D[h;  \tau \s]$ which is a  skew polynomial ring, $\s$ and $\tau$ are ring endomorphisms of $\CD$ that are extensions of the ring endomorphisms $\s$ and $\tau$ of $D$, respectively, defined  by the rule $\s (h) = \rho h +b$ and $\tau (h) = \tau (\rho^{-1}) (h-\tau (b))$. In particular, $\tau \s (h)=h$ and $\s\tau (h) = \o_\rho (h) = \rho \tau \s (\rho^{-1})h$. Furthermore, $\s\tau = \o_\rho \tau \s$ in $\CD$.
\end{theorem}
Theorem \ref{2Apr} is a generalization of a similar result for rings $D\langle \s ; b, \rho \rangle $, {\cite[Lemma 1.2, Corollary 1.4]{Bav-GlGWA-1996}}.

{\bf Simplicity criterion for diskew polynomial rings when $\rho$ is a unit.}
If $\rho$ is a unit, a simplicity criterion Theorem \ref{Di31Oct} for diskew polynomial rings is a relatively easy corollary of the simplicity criterion for generalized Weyl algebras (Theorem \ref{B5Apr}) because of Theorem \ref{2Apr}.
 \begin{theorem}\label{Di31Oct}%\marginpar{Di31Oct}
Let $E=D[x,y;\s , \tau , b, \rho ]$ be a diskew polynomial ring such that $\rho$ is a unit in $D$ and $\nu := \tau \s$ is an epimorphism. The following statements are equivalent.
\begin{enumerate}
\item The ring $E$ is a simple ring.
\item
\begin{enumerate}
\item The endomorphisms $\s$ and $\tau$ of $D$  are automorphisms,
\item the ring $D$ is a $\s$-simple ring,
 \item  for each natural number $n\geq 1$ there is no element $p=h^n+\sum_{i=0}^{n-1}d_ih^i\in \CD$,  where $d_i\in D$,  such that
\begin{enumerate}
\item for all elements $d\in D$, $pd= \nu^n (d) p$, i.e., $d_id= \nu^{n-i}(d) d_i$ for $i=0,1,\ldots , n-1$,
\item $\s (p) = \rho_n^\nu p$ where $ \rho_n^\nu =\rho \nu(\rho) \cdots \nu^{n-1}(\rho )$, and
\item $[h,p]= 0$, i.e, $\nu (d_i) =d_i$ for   $i=0,1,\ldots , n-1$, and
\end{enumerate}
\item the elements $b_i\in D$ (see (\ref{sih1})), where $i\geq 1$,   are units in $D$. In particular, $b=b_1\in D$ is a unit.
\end{enumerate}
\item
\begin{enumerate}
\item The endomorphisms $\s$ and $\tau$ of $D$ are automorphisms,
\item the ring $D$ is a $\tau$-simple ring,
 \item for each number $n\geq 1$ there is no element $p'=h'^n+\sum_{i=0}^{n-1}d_i'h'^i\in \CD =D[h', \mu := \s \tau ]$, where $d_i'\in D$ and $h'=\s (h)$, such that
\begin{enumerate}
\item for all elements $d\in D$, $p'd= \mu^n (d) p'$, i.e., $d_i'd= \mu^{n-i}(d) d_i'$ for $i=0,1,\ldots , n-1$,
\item $\tau (p') = (\rho^{-1})_n^\mu p$ where $ (\rho^{-1})_n^\mu:=\rho^{-1} \mu(\rho^{-1}) \cdots \mu^{n-1}(\rho^{-1} )$,     and
\item $[h',p']= 0$, i.e, $\mu (d_i') =d_i$ for   $i=0,1,\ldots , n-1$, and
\end{enumerate}
 \item the elements $b_i'\in D$ (see (\ref{sih2})), where $i\geq 1$,   are units in $D$. In particular, $b=b_1\in D$ is a unit.
\end{enumerate}
\end{enumerate}
\end{theorem}

{\it Remarks.} 1.  The conditions (a), (i) and (iii) in statement 2 imply that the element $p$ is a normal, regular element of $\CD$. So,  the ideal $(p) = p\CD = \CD p$ is a proper, $\s$-invariant ideal of the ring $\CD$ generated by the normal, regular element $p$ of $\CD$. The condition (ii) implies that $xp=\s (p) x= \rho^\nu_n px$ and $yp=\tau  (\rho^\nu_n)^{-1} py$, see  (\ref{pbbp1}). So, {\em the element $p$ is a normal, regular element of the ring $E$}.

2. In statement 2, the condition  (ii)  can be written as explicit equations on unknowns $d_0, d_1, \ldots , d_{n-1}$, see (\ref{hpdn7}).

3. Similar remarks can be made for the element $p'$ in statement 3 (by using the $(x,y)$-symmetry of GWAs, see Section \ref{GWA2}).

 Every simple ring is necessarily an algebra. Theorem \ref{QDi31Oct} and Theorem \ref{PDi31Oct} are refined versions of Theorem \ref{Di31Oct} in zero and prime characteristic, respectively.

{\bf Simplicity criterion for DPRs in characteristic zero.} If the ring $D$ is a $\Q$-algebra the condition (c) in Theorem \ref{Di31Oct} can be replaced by condition 4 of Theorem \ref{B6Apr}.
\begin{theorem}\label{QDi31Oct}%\marginpar{QDi31Oct}
Let $E=D[x,y;\s , \tau , b, \rho ]$ be a diskew polynomial ring such that $\rho$ is a unit in $D$,  $\nu = \tau \s$ is an epimorphism  and $D$ is a $\Q$-algebra. The following statements are equivalent.
\begin{enumerate}
\item The ring $E$ is a simple ring.
\item
\begin{enumerate}
\item The endomorphisms $\s$ and $\tau$ of $D$  are automorphisms,
\item the ring $D$ is a $\s$-simple ring,
 \item  there is no element $\alpha \in D$  such that $\rho \alpha - \s (\alpha ) = b$  and $\alpha d= \nu (d) \alpha$ for all elements $d\in D$.
 \item the elements $b_i\in D$ (see (\ref{sih1})), where $i\geq 1$,   are units in $D$. In particular, $b=b_1\in D$ is a unit.
\end{enumerate}
\end{enumerate}
\end{theorem}

{\bf Simplicity criterion for DPRs in prime  characteristic $p$.}
If the ring $D$ is a $\mF_p$-algebra the condition (c) in Theorem \ref{Di31Oct} can be replaced by more explicit conditions (where $\mF_p= \Z / p\Z$).
\begin{theorem}\label{PDi31Oct}%\marginpar{PDi31Oct}
Let $E=D[x,y;\s , \tau , b, \rho ]$ be a diskew polynomial ring such that $\rho$ is a unit in $D$,  $\nu = \tau \s$ is an epimorphism  and $D$ is a $\mF_p$-algebra. The following statements are equivalent.
\begin{enumerate}
\item The ring $E$ is a simple ring.
\item
\begin{enumerate}
\item The endomorphisms $\s$ and $\tau$ of $D$  are automorphisms,
\item the ring $D$ is a $\s$-simple ring,
 \item for each natural number $n\geq 0$ there is no element $p'=h^{p^n}+\sum_{i=0}^{n-1}\alpha_ih^{p^i}+\alpha$, where $\alpha , \alpha_i\in D$, such that
     \begin{enumerate}
\item for all $d\in D$, $pd=\nu^{p^n}(d)p$, i.e. $\alpha d = \nu^{p^n}(d) \alpha$ and $\alpha_i d = \nu^{p^n-p^i}(d) \alpha_i$ for $i=0,1,\ldots , n-1$,
\item $\s (p') =\rho^\nu_{p^n}p'$,
\item $[h,p']=0$, i.e. $\nu (\alpha ) = \alpha$ and $\nu (\alpha_i ) = \alpha_i$ for $i=0,1,\ldots , n-1$.
\end{enumerate}
\item the elements $b_i\in D$ (see (\ref{sih1})), where $i\geq 1$,   are units in $D$. In particular, $b=b_1\in D$ is a unit.
\end{enumerate}
\item
\begin{enumerate}
\item The endomorphisms $\s$ and $\tau$ of $D$  are automorphisms,
\item the ring $D$ is a $\s$-simple ring,
 \item  there is no element $\alpha \in D$  such that $\rho \alpha - \s (\alpha ) = b$, $\nu (\alpha ) = \alpha $  and $\alpha d= \nu (d) \alpha$ for all elements $d\in D$, and for each natural number $n\geq 1$ there are no elements $\alpha, \alpha_0, \ldots , \alpha_n$ such that
     \begin{enumerate}
\item for all $d\in D$, $\alpha d = \nu^{p^n}(d) \alpha$ and $\alpha_i d = \nu^{p^n-p^i}(d) \alpha_i$ for $i=0,1,\ldots , n-1$,
\item $\s (\alpha_i) =\rho^\nu_{p^n-p^i}\alpha_i$ for $i=0,1,\ldots , n-1$, and  $\rho^\nu_{p^n}\alpha -\s (\alpha ) = b^{p^n}+\sum_{i=0}^{n-1}\s(\alpha_i) b^{p^i}$,
\item $\nu (\alpha ) = \alpha$ and $\nu (\alpha_i ) = \alpha_i$ for $i=0,1,\ldots , n-1$.
\end{enumerate}
\item the elements $b_i\in D$ (see (\ref{sih1})), where $i\geq 1$,   are units in $D$. In particular, $b=b_1\in D$ is a unit.
\end{enumerate}
\end{enumerate}
\end{theorem}
{\bf The canonical left normal element $C$ of a diskew polynomial ring}.
 Theorem \ref{B6Apr} is a criterion for an element $C=h+\alpha$ (where $\alpha \in D$) to be a left normal element in $E$, it is a key moment in the proof of Theorem \ref{QDi31Oct} and Theorem \ref{PDi31Oct} (together with ``Meeting the $p$-neighbour method', see the proof of Theorem \ref{PDi31Oct}). Theorem \ref{B6Apr} is a generalization of {\cite[Lemma 1.3]{Bav-GlGWA-1996}}.

\begin{theorem}\label{B6Apr}%\marginpar{B6Apr}
Let $E= D[x,y; \s , \tau , b, \rho ]$ be a diskew polynomial ring such that $\rho$ is a unit, $\CD = D[h; \nu = \tau \s ]$ and $C= h+\alpha$ where $h=yx$ and $\alpha \in D$.  The following statements are equivalent.
\begin{enumerate}
\item The element $C$ is left normal in $E$.
\item $\rho \alpha - \s (\alpha ) = b$, $\nu (\alpha )=\alpha $  and $\alpha d= \nu (d) \alpha$ for all elements $d\in D$.

If one of the equivalent conditions holds then $[h,C]=0$ and
\begin{enumerate}
\item  $C=\rho^{-1} (xy+\s (\alpha ))$, $xC= \rho Cx$ and $yC= \tau (\rho^{-1}) C y$.
\item $E\simeq D[C; \nu ] [x,y; \s , \tau , a:= C-\alpha ]$ is a GWA where $\s (C) = \rho  C$ and $\tau (C) = \tau (\rho^{-1}) C$.
\item The element $C$ is a left normal, left regular element of $E$ and
 $E/ (C)\simeq D [x,y; \s , \tau , -\alpha ]$ is a GWA.
\item The element $C$ is a normal element in $E$ iff $\im (\nu ) = D$.
\item The element $C$ is regular  iff  $C$ is  right regular iff $\ker (\nu  ) =0$.
\item The element $C$ is a normal,  regular element iff $\nu$ is an automorphism of $D$.
\end{enumerate}
\end{enumerate}
\end{theorem}

%%%%%%%%%%%%%%%%%% SECTION 2 %%%%%%%%%%%%%%%%%%%%%%%%

\section{Generalized Weyl algebras}\label{GWA2}%\marginpar{GWA2}

At the beginning of the section  we consider examples of generalized Weyl algebras. We show that each GWA $A$ is a $\Z$-graded ring which is a free left $D$-module (Theorem \ref{1Apr})  but not a free right $D$-module, in general. Certain important left denominator sets of $A$ are considered in Proposition \ref{b1Apr} and Proposition \ref{d2Apr}. Results of this section are used throughout the paper.

{\bf Left and right normal elements}. Let $D$ be a ring. An element $a\in D$ is called a {\em left} (resp., {\em right}) {\em normal} element of $D$ if $aD\subseteq Da$ (resp., $Da\subseteq aD$).  If $a$ is a left normal element in $D$ then the left ideal $Da$ is an ideal of $D$. Similarly, if $a$ is a right normal element in $D$ then the right ideal $aD$ is an ideal of $D$.  An element $a\in D$ is {\em normal}  if $aD =Da$, i.e., $a$ is left and right normal. Let $\cdot a:=\cdot a_D: D\ra D$, $d\mapsto da$,   and  $\ga := \ker (\cdot a)$. In particular, $\ga a=0$.  Similarly, let $a\cdot :=a_D\cdot : D\ra D$, $d\mapsto ad$,   and  $\gb := \ker (a\cdot )$. In particular, $a\gb =0$.
 If the element $a$ is left normal then $\gb$ is an ideal of the ring $D$: $a\cdot D\gb D\subseteq Da\gb D=0$.   If the element $a$ is right  normal then $\ga$ is an ideal of the ring $D$: $D\ga D\cdot a\subseteq D\ga a D=0$. The sets $\mL_a:=\{ d\in D\, | \, da=ad'$ for some $d'\in D\}$ and $\mR_a:=\{ d\in D\, | \, ad=d'a$ for some $d'\in D\}$ are subrings of $D$ such that $\ga \subseteq \mL_a$ and $\gb \subseteq \mR_a$. Furthermore, $\ga$ is an ideal of $\mL_a$ ($\mL_a\ga \, \mL_a\cdot  a\subseteq \mL_a\ga\,  aD=0$, and so $\mL_a\ga \, \mL_a\subseteq \ga$) and $\gb$ is an ideal of $\mR_a$ ($a\cdot \mR_a\gb \mR_a\subseteq D a\gb \mR_a=0$). If $a$ is a left (resp., right) normal element of $D$ then $\mL_aa=aD$ (resp., $Da= a\mR_a$).

Suppose that $a\in D$ is a left normal element. Then, for each element $d\in D$, $ad= d_la$ for some element $d_l\in \mL_a$ which is unique up to adding $\ga$ ($d_la= (d_l+\ga ) a$). Hence, the map
%\marginpar{waLa}
\begin{equation}\label{waLa}
\o_a: D/\gb \ra \mL_a/ \ga, \;\; d+\gb \mapsto d_l+\ga,
\end{equation}
is a ring {\em isomorphism}  and we can write $ad=\o_a(d)d$ for all $d\in D$. A left normal element $a$ is normal iff $\mL_a=D$. If  $a$ is a normal element  then $\mL_a=D$ and the map
$\o_a: D/\gb \ra D/ \ga$, $ d+\gb \mapsto d_l+\ga$,
is a ring isomorphism.

Similarly, suppose that $a\in D$ is a right normal element. Then, for each element $d\in D$, $da= ad_r$ for some element $d_r\in \mR_a$ which is unique up to adding $\gb$ ($ad_r=a (d_r+\gb ) $). Hence, the map
%\marginpar{waRa}
\begin{equation}\label{waRa}
\o_a': D/\ga \ra \mR_a/ \gb, \;\; d+\ga \mapsto d_r+\gb,
\end{equation}
is a ring {\em isomorphism} and we can write $da= a\o_a'(d) $ for all $d\in D$.

 A right normal element $a$ is  normal iff $\mR_a=D$. If  $a$ is a normal element then $\mR_a=D$ and the map
$\o_a': R/\ga \ra R/ \gb$, $ d+\ga \mapsto d_r+\gb$,
is a ring isomorphism.  So, for a normal element $a$ of $D$,
%\marginpar{adwda}
\begin{equation}\label{adwda}
ad=\o_a(d)a\;\; {\rm and}\;\; da=a\o_a'(d)\;\; {\rm   for \; all}\;\; d\in D.
\end{equation}

\begin{lemma}\label{a2Apr}%\marginpar{a2Apr}
If $a\in D$ is a normal element then the maps $\o_a: D/ \gb \ra D/ \ga$ and $\o_a': D/ \ga \ra D/ \gb$ are ring isomorphisms such that $\o_a'= \o_a^{-1}$. If, in addition, $a$ is a regular element then the maps $\o_a, \o_a': D\ra D$ are ring isomorphisms such that $\o_a'=\o_a^{-1}$.
\end{lemma}

{\it Proof}. The lemma follows from (\ref{adwda}): For all $d\in D$, $ad= \o_a(d) a =a \o_a'\o_a(d) $ and $ da= a\o_a'(d) = \o_a\o_a'(d) a$.  $\Box $

{\bf Examples of generalized Weyl algebras where $a$ is central, \cite{Bav-GWA-FA-91}-\cite{Bav-Bav-Ideals-II-93}.} 1. The (first) {\em Weyl algebra} $A_1= K\langle x, \der \, | \, \der x -x\der =1\rangle $ over a ring $K$ is the GWA $K[h][x,y:=\der ; \s , a=h]$ with base ring $K[h]$ and its $K$-automorphism $\s$ defined by the rule $\s (h) = h-1$.

2. The {\em quantum plane} $\L = K\langle x,y \, | \, xy= qyx\rangle$ where $q$ is a central unit of $K$ is the GWA $K[h][x,y; \s , a=h]$ where $\s (h) = qh$.

3.  For $q,\,h=q-q^{-1}\in K=\mathbb{C} $, the algebra
$U_q=U_qsl(2)$ is generated by $X, Y, H_-$ and  $H_+$ subject to the defining relations:
$$H_+H_-=H_-H_+=1,\;\; XH_{\pm}=q^{\pm 1}H_{\pm}X,\;\; YH_{\pm}=
q^{\mp 1}H_{\pm}Y, \;\; [X,Y]=\frac{H_+^2-H_-^2}{h}.$$
It follows that the algebra $U_q$  is a GWA,
$$U_q\simeq K[C,H,H^{-1}](\sigma,a=C+\Big(H^2/(q^2-1)-H^{-2}/(q^{-2}-1)\Big)/2h),$$
where $\sigma(H)=qH,\  \sigma(C)=C$.

4. {\em Woronowicz's deformation} $V$ is an algebra generated by elements $V_0$, $V_-$ and  $V_+$ subject to the defining relations,  \cite{Za}:
$$[V_0,V_+]_{s^2}:= s^2V_0V_+-s^{-2}V_+V_0=V_+,
\,\,[V_-,V_0]_{s^2}=V_-,\;\;
[V_+,V_-]_{1/s}:=  s^{-1}V_+V_- -sV_-V_+=V_0.$$
 The algebra  $V$ is a GWA,
$$V\simeq K[u,v](\s , a=v),\;\; V_+\lra x, \;\;V_-\lra y,\;\;V_0\lra u,\;\;
V_-V_+\lra v,$$
where $\s :u\ra s^2(s^2u-1)$, $v\ra s^2v+su$, is the  automorphism of the
polynomial ring $K[u,v]$ in two variables $u$ and $v$. Let
$$H=u+s^2/(1-s^4)\;\; {\rm and}\;\; Z=v+u/s(1-s^2)+s^3/(1-s^2)(1-s^4).$$
Then $\s (H)=s^4H$, $\s (Z)=s^2Z$ and $K[u,v]=K[H,Z]$. So,
$$V\simeq K[H,Z](\s,a=Z+\alpha H+\b ),\;\; V_+\lra x, \;\;V_-\lra y,\;\;
 V_0\lra H-s^2/(1-s^4),$$
where $\s :H\ra s^4H$, $Z\ra s^2Z$; $\alpha =-1/s(1-s^2)$ and $\b =s/(1-s^4)$.

5. {\em Witten's first deformation}  $E$ is an  algebra generated by elements
$E_0$, $E_-$ and  $E_+$  subject to the defining relations, \cite{Za}:
$$[E_0,E_+]_p:= pE_0E_+-p^{-1}E_+E_0=E_+,
\;\; [E_-,E_0]_p=E_-,\;\;
[E_+,E_-]=E_0-(p-1/p)E^2_0,$$
where $p\not= 0,\pm 1,\pm i\in K$. The element $C=E_-E_++ \frac{E_0(E_0+p)}{p(p^2+1)}$ is central in $E$.
  Witten's first deformation is a GWA,
 $$E\simeq K[C,H](\s, a=C-H(H+1)/(p+p^{-1})),\;\; E_+\lra x, \;\;E_-\lra y,\;\;
  E_0\lra pH,$$
where $\s :C\ra C$, $H\ra p^2(H-1)$.

6. {\em The quantum group } ${\cal O}_{q^2}(so(K,3))=K[H]\langle \s ;
b=(q-q^{-1})H,\rho =1\rangle ,$ $\s (H)=q^2H$, $q\in K$, \cite{Smith1}, %by Corollary 1.6
is isomorphic to the GWA of degree 1:
$${\cal O}_{q^2}(so(k,3))=K[H,C](\s ,a=C+H^2/q(1+q^2)),\,\s (H)=q^2H,\,
\s (C)=C.$$

{\bf Generalized Weyl algebras with 2 endomorphisms and a left normal element $a$.} Let $A= D[x,y; \s , \tau , a]$ be a generalized Weyl algebra.

{\em  Remarks}. 1. If either $\s$ or $\tau$ is an automorphism then the last two conditions in (\ref{abGWA})  are equivalent provided $\tau \s (a) =a$. In more detail, if $\s$ is an automorphism then the third condition in (\ref{abGWA})  is obtained from the second one by applying $\s$ and replacing $d$ by $\s (d)$. If $\tau$ is an automorphism then the second condition in (\ref{abGWA}) is obtained from the third by applying $\tau$, using the equality $\tau \s (a) =a$ and replacing $\tau (d)$ by $d$.

2. If $\s$ and $\tau$ are automorphisms and $\tau = \s^{-1}$ then (\ref{abGWA}) is equivalent to the fact that the element $a$ belongs to the centre of $D$ and we obtain the usual definition of GWA.

3. Let $a$ be a normal element of $D$ ($da=aD$) which is also a regular element of $D$ ($a$ is a non-zero-divisor). Then for all elements $d\in D$, $ad=\o_a(d)a$ for an automorphism $\o_a$ of $D$. Let $\tau = \o_a\s^{-1}$. Then $A= D[x,y; \s ,\tau , a] $ is a GWA.

{\bf Examples of generalized Weyl algebras where $a$ is left normal.} 1. Let $K$ be a field and $\l \in K^*:= K\backslash \{ 0\}$. Let $D$ be the {\em quantum plane} $\L := K\langle
p,q\, | \, pq=\l qp\rangle$.  If $\l$ is not
a root of unity then the set of all nonzero normal elements in
$\L$ is $\{ K^* p^iq^j\, | \, i,j\geq 0\}$ and the centre $Z$ of
$\L$ is $K$. If $\l$ is  a primitive $n$'th root of unity ($n\geq
2)$ then the centre $Z$ of $\L$ is the polynomial algebra $K[p^n,
q^n]$ in two variables $p^n$ and $q^n$, and
 the set of all nonzero normal elements in $\L$ is $\{ Z^*
p^iq^j\, | \, 0\leq i,j<n\}$ where $Z^*:= Z\backslash \{ 0\}$.
Each pair $ (\alpha , \beta )\in K^2$ of nonzero scalars
determines the $K$-automorphism $\s =\s_{\alpha , \beta }$ of the
algebra $\L$: $p\mapsto \alpha p$, $ q\mapsto \beta q$. Let $ a=
zp^iq^j$ where $z\in Z^*$ and $i,j\geq 0$, if $\l $ is not  a root of
unity;  or $0\leq i,j<n$, if $\l $ is an $n$'th primitive root of
unity. The element $a$ is a normal regular element and $\o_a (p)=
\l^{-j}p$, $ \o_a(q)= \l^iq$. Then $\o_a\s^{-1}(p)= \alpha^{-1}
\l^{-j} p$ and $\o_a\s^{-1} (q) = \beta^{-1} \l^iq$.  Then the GWA
$A= \L [x,y;\s , \tau = \o_a\s^{-1}, a]$ is a $K$-algebra generated by $\L$, $x$ and $y$
that are subject to the defining relations:
\begin{eqnarray*}
xp=\alpha px, \;\;\;   xq=\beta qx,\; &\;  yp= \alpha^{-1} \l^{-j} py,\;\;\;  yq= \beta^{-1} \l^i qy,  \\
 yx= zp^iq^j, & xy= \alpha^i \beta^j \s (z) p^iq^j.
\end{eqnarray*}
The algebra $A$ is a domain of Gelfand-Kirillov dimension $3$.

2. Let $D=C[t_1, \ldots , t_n; \nu_1, \ldots , \nu_n]$ be a skew polynomial ring in variables $t_1, \ldots , t_n$ ($t_it_j=t_jt_i$, $t_ic=\nu_i(c)t_i$ for all $c\in C$) over a ring $C$, $\nu_1, \ldots , \nu_n$ are commuting automorphisms of the ring $C$. Suppose that an element $u\in C$ is such that $\nu_i(u) = u_iu$ for some unit $u_i$ of $D$ for $i=1, \ldots , n$ (eg, $u=u_1=\cdots =u_n=1$). Then for each element $\alpha = (\alpha_1, \ldots , \alpha_n)\in \N^n$, the element $a=ut^\alpha$ is a regular, normal element of $D$ where $t^\alpha =t_1^{\alpha_1}\cdots t_n^{\alpha_n}$. Then $A= D[x,y;\s, \o_a\s^{-1}, a]$ is a GWA.

{\bf The right generalized Weyl algebras.} Let $A$ be a ring. The {\em opposite ring} $A^{op}$ of $A$ is a ring that is equal to $A$ as an abelian group but the multiplication in $A^{op}$ is given by the rule $a\cdot b = ba$. The definition of GWA is not left-right symmetric. By taking the opposite ring of a GWA we obtain the definition of a right generalized Weyl algebra.

{\it Definition}. Let $D$ be a ring, $\s $ and $\tau$ be ring endomorphisms of the ring $D$, and an element $a\in D$ be  such that
%\marginpar{rabGWA}
\begin{equation}\label{rabGWA}
\tau \s (a) = a, \;\; da= a\tau \s (d) \;\; {\rm and}\;\; d\s (a) =\s (a) \s\tau (d)  \;\; {\rm for \; all}\;\; d\in D.
\end{equation}
The {\bf right generalized Weyl algebra} (rGWA) $A'= D(\s, \tau, a)_r = D[x,y; \s, \tau, a]_r$ is a ring generated by $D$, $x$ and $y$ subject to the defining relations:
%\marginpar{rGWADEF}
\begin{equation}\label{rGWADEF}
dx=x\s (d) \;\;  {\rm and} \;\;  dy=y\tau (d)\;\; {\rm for \; all} \;\; d\in D,
 \;\; xy=a \;\; {\rm and} \;\;   yx=\s (a).
\end{equation}
The ring $D$ is called the {\em base ring} of the rGWA $A$. The endomorphisms $\s$, $\tau$ and  the element $a$ are called the {\em
defining endomorphisms} and the  {\em defining element} of the rGWA $A$, respectively. By (\ref{rabGWA}), the elements $a$ and $\s (a)$ are right normal in $D$.

{\it Example.} Let $A= D[x,y; \s, \tau, a]$ be a GWA. Then $A^{op}= D^{op}[x,y; \s, \tau, a]_r$
 is a rGWA. Similarly, if $A'= D[x,y; \s, \tau, a]_r$ is a rGWA. Then $A^{op}= D^{op}[x,y; \s, \tau, a]$ is a GWA. In this paper, we study GWAs. Analogous properties of rGWAs are obtained by applying the functor $A\ra A^{op}$. In the obvious way, iterated rGWAs are defined.

{\bf A $\Z$-grading of a GWA.} The next theorem proves existence of GWAs and introduces a $\Z$-grading.

\begin{theorem}\label{1Apr}%\marginpar{1Apr}
The GWA $A= D[x,y; \s ,\tau , a]$ exists. It is a $\Z$-graded ring $A= \oplus_{i\in \Z} A_iv_i$ where $A_i = Dv_i\simeq {}_DD$,  $v_0=1$, $v_i = x^i$ and $v_{-i} = y^i$ for $i\geq 1$. In particular, the module ${}_DA$ is  free.
\end{theorem}

{\it Proof}. Consider a free left $D$-module $A' := \oplus_{i\in \Z} Dv_i'$ (it is a prototype of $A$) where elements  $v_i'$ are free left $D$-module generators. The aim is to define an action of the elements $x$ and $y$ on $A'$  such that the relations (\ref{GWADEF}) hold and therefore to realize the ring $A$ as a subring of $\End_\Z (A')$. This proves existence of the ring $A$. Then using the $A$-module $A'$ we show that the ring $A$ is a $\Z$-graded ring (in fact, ${}_AA\simeq {}_AA'$).

The action is given by the rule: For $d'\in D$ and $i\in \Z$,
$$ x\cdot d'v_i':= \begin{cases}
\s (d') v_{i+1}' & \text{if }i\geq 0,\\
 \s (d') \s (a) v_{i+1}' & \text{if }i<0, \\
\end{cases} \;\;\;\; \;\;\;
 y\cdot d'v_i':= \begin{cases}
\tau (d') v_{i-1}' & \text{if }i\leq 0,\\
\tau (d') a v_{i-1}' &\text{if }i>0. \\
\end{cases}
$$
Let us check that the relations (\ref{GWADEF}) hold.

(i) $xd=\s (d) x$: $xd\cdot d'v_i'= \s (d) \s (d') \cdot \begin{cases}
v_{i+1}' & \text{if }i\geq 0,\\
\s (a) v_{i+1}' & \text{if }i<0, \\
\end{cases}= \s (d) x\cdot d'v_i'.$

(ii) $yd=\tau (d) x$: $yd\cdot d'v_i'= \tau (d) \tau (d') \cdot \begin{cases}
v_{i-1}' & \text{if }i\leq 0,\\
a v_{i-1}' & \text{if }i>0, \\
\end{cases}= \tau (d) y\cdot d'v_i'.$

(iii) $yx=a$: $yx\cdot d'v_i'= y  \cdot \s (d') \cdot \begin{cases}
v_{i+1}' & \text{if }i\geq 0,\\
\s (a) v_{i+1}' & \text{if }i<0, \\
\end{cases}= \tau \s (d') \cdot \begin{cases}
av_i' & \text{if }i\geq 0,\\
\tau \s (a) v_i' & \text{if }i<0, \\
\end{cases}=\tau \s (d') av_i'= a\cdot d'v_i'$  (since $\tau \s (a)=a$ and $ad' = \tau\s (d') a$).

(iv) $xy=\s (a)$: $xy\cdot d'v_i'= x  \cdot \tau (d') \cdot \begin{cases}
v_{i-1}' & \text{if }i\leq 0,\\
a v_{i-1}' & \text{if }i>0, \\
\end{cases}= \s\tau  (d') \cdot \begin{cases}
\s (a)v_i' & \text{if }i\leq 0,\\
\s (a) v_i' & \text{if }i>0, \\
\end{cases}=\s\tau (d') \s (a)v_i'= \s (a)\cdot d'v_i'$ (since $\s (a)d' = \s\tau (d') \s (a)$). So, the ring $A$ exists.

By (\ref{GWADEF}), the ring  $A$ is equal to the sum $\sum_{i\in \Z} A_i$ where $A_i = Dv_i$. Since $A_i \cdot v_0'= Dv_i'\simeq {}_DD$, we have $A\cdot v_0'=\oplus_{i\in \Z}A_i v_0'$. Therefore, $A= \oplus_{i\in \Z}A_i$ and $A_i = Dv_i\simeq {}_DD$ for all $i\in \Z$. So, the left $D$-module  $A$ is free. $\Box $

For all $n,m\in \Z$, %\marginpar{1vnvmm}
\begin{equation}\label{1vnvmm}
v_nv_m= (n,m) v_{n+m}
\end{equation}
for some elements $(n,m)\in D$ where, for all $n>0$ and $m>0$,
\begin{eqnarray*}
 n\geq m: & (n, -m)= \s^n(a) \cdots \s^{n-m+1} (a), & (-n, m)= \tau^{n-1}(a) \cdots \tau^{n-m} (a), \\
 n\leq m: & (n, -m)= \s^n(a) \cdots \s (a), & (-n, m)= \tau^{n-1}(a) \cdots
 \tau (a)a.
\end{eqnarray*}
For all other values of $n$ and $m$, $(n,m)=1$.
By Theorem \ref{1Apr}, each element $r$  of the ring $A$ is a unique (finite) sum $r=\sum_{i\in \Z} r_iv_i$ where $r_i\in D$. If the element $r$ is  nonzero  then the natural number $l(r):=n-m$ is called the {\em length} of $r$ where $n=\max \{ i\in \Z\, | \, r_i\neq 0\}$ and  $m=\min \{ i\in \Z\, | \, r_i\neq 0\}$. The multiplication in the ring $A$ is given by the rule: for all $d,d'\in D$ and $i,j\in \Z$,
%\marginpar{dmult1}
\begin{equation}\label{dmult1}
dv_i\cdot d'v_j = \begin{cases}
d\s^i (d') (i,j) v_{i+j} & \text{if }i\geq 0,\\
d\tau^i(d') (i,j) v_{i+j} & \text{if }i<0. \\
\end{cases}
\end{equation}
For all $i\geq 1$,
%\marginpar{iist}
\begin{equation}\label{iist}
(i,-i) = \s^i((-i,i))\;\; {\rm and}\;\; (-i,i) = \tau^i((i,-i)).
\end{equation}
{\it Proof}. The equalities $(i,-i) x^i = x^iy^ix^i= x^i (-i,i) = \s^i((-i,i))x^i$ imply that $(i,-i) = \s^i((-i,i))$ (by Theorem \ref{1Apr}). Similarly, the equalities $(-i,i) y^i = y^ix^iy^i= y^i (i,-i) = \tau^i((i,-i))y^i$ imply that
$(-i,i) = \tau^i((i,-i))$ (by Theorem \ref{1Apr}). $\Box$

Clearly, the equalities (\ref{iist}) hold in the ring $D$. So, they must follow from the equalities (\ref{abGWA}). It is not that straightforward to prove them in this way. For example, for $i=2$ we have $\tau^2((2,-2)) = (-2, 2)$, that is $\tau^2(\s^2(a)\s(a))=\tau (a) a$. To prove this equality using (\ref{abGWA}) we proceeds as follows:
$$ \tau^2(\s^2(a)\s(a))=\tau (\tau \s^2(a)\tau \s(a))\stackrel{(\ref{abGWA})}{=}\tau (\tau \s^2 (a) a) \stackrel{(\ref{abGWA})}{=} \tau (a\s (a))=\tau (a) \tau \s (a)
\stackrel{(\ref{abGWA})}{=}\tau (a) a.$$

For a natural number $i\geq 1$, let $\s_i:=\s^i$ and $\s_{-i}:=\tau^i$, and $\s_0:= \id_D$, the identity automorphism of $D$.  For all $m,n\in \Z$ and $d\in D$,
%\marginpar{1iist}
\begin{equation}\label{1iist}
\s_n\s_m(d) (n,m)=(n,m)\s_{n+m}(d).
\end{equation}
This follows from (\ref{1vnvmm}) and Theorem \ref{1Apr}. By (\ref{1iist}), {\em  the elements $(-i,i)$ and $(i,-i)$ are left normal in $D$} (since $\s_0= \id_D$).

{\bf The $(x,y)$-symmetry of a GWA}. By applying the endomorphism $\s$ to the equality $\tau \s (a) =a$  and then interchanging the last two equalities in (\ref{abGWA}), we obtain that (since $\tau \s (a)=a$)
$$\s\tau \s (a) = \s (a), \;\; \s (a) d= \s\tau (d) \s (a)\;\; {\rm and}\;\; \tau\s (a) d= \tau \s (d) \tau \s (a)\;\; {\rm for \; all}\;\; d\in D.$$
So, we have the GWA $D[x',y'; \tau , \s , \s (a)]$. Clearly,
%\marginpar{Asyma1}
\begin{equation}\label{Asyma1}
D[x,y; \s , \tau , a]=D[y,x; \tau , \s , \s (a)].
\end{equation}
The `identity' isomorphism ($o: u\mapsto u^o$)
%\marginpar{Asyma}
\begin{equation}\label{Asyma}
 o : D[x,y; \s , \tau , a]\ra D[x',y'; \tau , \s , \s (a)], \;\; x\mapsto y', \;\; y\mapsto x', \;\; d\mapsto d\;\; (d\in D),
\end{equation}
is called the (canonical) $(x,y)$-{\em symmetry isomorphism}  or the {\em symmetry isomorphism}, for short. Clearly, $r^{oo}=r$ for all $r\in A$.  The symmetry isomorphism is a $D$-isomorphism that  reverses the $\Z$-grading: $A_i^{o}=A_{-i}$ for all $i\in \Z$.

{\bf Ore and denominator sets}. Let $S$ be a nonempty subset of a ring $R$. Let $\ass_l (S):=\{ r\in R\, | \, sr=0$  for some $s=s(r) \in S\}$ and $\ass_r (S):=\{ r\in R\, | \, rs=0$  for some $s=s(r) \in S\}$. A nonempty subset $S$ of  $R\backslash \{ 0\}$ is called a {\em multiplicative set} if $SS\subseteq S$ and $1\in S$. A multiplicative set $S$ is called a {\em left Ore set} (resp., a {\em right Ore set}) if for given elements $s\in S$ and $r\in R$, $Sr\cap Rs\neq \emptyset$ (resp., $rS\cap sR\neq \emptyset $).
 If $S$ is a left (resp., right) Ore set of $R$ then $\ass_l(S)$ (resp., $\ass_r(S)$) is an ideal of the ring $R$. The sets of left and right Ore sets of $R$ are denoted by $\Ore_l(R)$ and $\Ore_r(R)$, respectively. Their intersection $\Ore (R) = \Ore_l(R)\cap \Ore_r(R)$ is the {\em set of Ore sets} of $R$.

 A left Ore set $S$ of $R$ is called a {\em left denominator set} of $R$ if $\ass_l(S)\supseteq \ass_r(S)$. Similarly, a right Ore set $S$ of $R$ is called a {\em right denominator set} of $R$ if $\ass_l(S)\subseteq \ass_r(S)$. The sets of left and right denominator sets are denoted $\Den_l(R)$ and $\Den_r(R)$, respectively. Their intersection $\Den (R) = \Den_l(R)\cap \Den_r(R)$ is the {\em set of denominator sets}. For an ideal $\ga$ of $R$, $\Den_l(R, \ga ) :=\{ S\in \Den_l(R)\, | \, \ass_l(R)=\ga \}$ and $\Den_r(R, \ga ) :=\{ S\in \Den_r(R)\, | \, \ass_r(R)=\ga \}$. For each $S\in \Den_l(R)$, the ring $S^{-1}R=\{ s^{-1}r\, | \, s\in S, r\in R\}$ is called the {\em left quotient ring} of $R$ at $S$ or the {\em left localization} of $R$ at $S$. For each $S\in \Den_r(R)$, the ring $RS^{-1}=\{ rs^{-1}\, | \, s\in S, r\in R\}$ is called the {\em right  quotient ring} of $R$ at $S$ or the {\em right localization} of $R$ at $S$. Let $\Den (R):= \Den_l (R)\cap \Den _r(R)$ and $S\in \Den (R)$. Then $\ass_l(S) = \ass_r(S)$ and $S^{-1}R\simeq RS^{-1}$.

{\bf The sets $S_x$ and $S_y$ are left denominator sets}.  Let $R$ be a ring and $\s : R\ra R$ be a ring endomorphism. Then $\ker (\s) \subseteq \ker (\s^2)\subseteq \cdots \subseteq \ker (\s^i ) \subseteq \cdots $ is an ascending chain of ideals of $R$. Their union  $\CK (\s ) := \bigcup_{i\geq 1} \ker (\s^i )$ is an ideal of $R$ such that $\s (\CK (\s )) \subseteq \CK (\s )$. Let $R(\s ) := R/ \CK (\s )$. The map
%\marginpar{RsKs}
\begin{equation}\label{RsKs}
\overline{\s} : R(\s ) \ra R(\s ) , \;\; r+\CK (\s ) \mapsto \s (r)+ \CK (\s ),
\end{equation}
is a ring {\em monomorphism}. Let $R$ be either a free algebra $K\langle x_0, x_1, \ldots \rangle$ or a polynomial algebra $K[x_0, x_1, \ldots ]$ in countably many variables over a ring $K$ and $\s$  an $K$-endomorphism of $R$ given by the rule
 $\s (x_0) =0$ and $\s (x_i) = x_{i-1}$ for all $i\geq 1$. Then $\CK (\s)$ is the ideal $(x):= (x_0, x_1, \ldots)$ of $R$, $R(\s ) \simeq K$ and $\overline{\s}: K\ra K$ is the identity map.

Let $I$ be a nonempty subset of $R$, $r\in R$ and $(I:r):=\{ u\in R\, | \, ur \in I\}$. If $I$ is a left ideal of the ring $R$ then so is $(I:r)$. If $M$ (resp., $N$) is a left (resp., right) $R$-module we also write ${}_RM$ (resp., $N_R$) to indicate this fact. For  $r\in R$, let $r_M\cdot :M\ra M$, $m\mapsto rm$, and $\cdot r_N: N\ra N$, $n\mapsto ur$. Let $D$ be a ring. An element $d\in D$ is called {\em left} (resp., {\em right} ) {\em regular} if $\ker (\cdot d_D)=0$ (resp., $\ker ( d_D\cdot )=0$). A left and right regular element is called {\em regular}. Let $\pCC_D$, $\CC_D'$ and $\CC_D$ be the sets of left, right and regular elements of $D$, respectively. For  a subset $C$ of $D$, the set $\s^{-i}(C):= \{ d\in D\, | \, \s^i (d) \in C\}$ is the pre-image of the set $C$ under the map $\s^i : D\ra D$. For simplicity reason, we use the notation $\s^{-i}(C)$ rather than $(\s^i)^{-1}(C)=\underbrace{\s^{-1}\cdots \s^{-1}}_{i \, {\rm times}} (C)$ to denote the pre-image of $C$.

\begin{proposition}\label{b1Apr}%\marginpar{b1Apr}
Let $A= D[x,y; \s , \tau , a]$ be a GWA and $S_x=\{ x^i \, | \, i\geq 0\}$. Then
\begin{enumerate}
\item $S_x\in \Den_l(A, \ga )$ where $\ga := \bigoplus_{i\geq 1} \s^{-i} (\CK  (\s) :(i,-i))y^i \oplus \bigoplus_{i\geq 0} \CK (\s ) x^i$ and the set \;\;\;\;\;\;  $\s^{-i}(\CK (\s ) :(i,-i))$ is equal to $ \{ d\in D\, | \, \s^i(d)(i,-i)\in \CK (\s ) \}$.
\item $A/ \ga \simeq \bigoplus_{i\geq 1} D/( \s^{-i} (\CK  (\s) :(i,-i)))y^i \oplus \bigoplus_{i\geq 0} D (\s ) x^i$ where $D(\s ) := D/\CK (\s )$.
\item $\ass_r(S_x)=\bigoplus_{i\geq 1} \ker (\cdot (-i,i)_D)y^i\subseteq \ga$.
\item $A_x:= S_x^{-1}A=A_{x,0}[x^{\pm 1}; \s ]$ is a skew polynomial ring where $A_{x,0}:= \bigcup_{i\geq 0} x^{-i}D (\s )x^i$ and $\s (x^{-i} (d+\CK (\s )) x^i) := x^{-i} (\s (d)+\CK (\s )) x^i$. The addition and multiplication in the ring $A_{x,0}$ are given in (\ref{xidxi1}).
    \item
     \begin{enumerate}
\item $\ker_A(x\cdot ) = \bigoplus_{i\geq 1}(\ker (\s) : a)y^i\oplus\bigoplus_{i\geq 0}\ker (\s ) x^i$; $x\in \CC_A'$ iff $\s$ a monomorphism  and $a\in \lCCD $.
\item $\ker_A(\cdot x) = \bigoplus_{i\geq 1}(\ker_D (\cdot \tau^{i-1}(a))y^i$; $x\in \lCCA$ iff $\tau^i(a)\in \lCCD$ for all $i\geq 0$.
\item $x\in \CC_A$ iff $x\in \CC_A'$ iff $\s$ is a monomorphism and $a\in \lCCD$.
\end{enumerate}
\item $\ass_l(S_x)=\ass_r(S_x)$ iff $\s$ is a monomorphism.
\item $S_x\in \Den (R)$ iff $\s$ is an automorphism. If $S_x\in \Den (R)$ then $\ass (S_x) =\bigoplus_{i\geq 1} \ker (\cdot (-i,i)_D)y^i$.
\end{enumerate}
\end{proposition}

{\it Proof}. 3. The element $x$ is a homogeneous element of the $\Z$-graded ring $A= \oplus_{i\in \Z} A_i$. Therefore, $\gb := \ass_r(S_x) = \oplus_{i\in \Z} \gb_i$ where $\gb_i := \gb \cap A_i$. By Theorem \ref{1Apr}, for $i\geq 0$, $\gb_i=0$, since for all $j\geq 1$ the map $\cdot x^j: Dx^i\ra Dx^{i+j}$ is an injection. Now, $\gb_{-i} = \ker (\cdot (-i,i)_D)$ since $ \cdot x^{i+j} : A_{-i} =Dy^i\ra Dx^j$, $dy^i\mapsto d(-i,i) x^j$. The inclusion $\gb \subseteq \ga$ is proven in the proof of statement 1 (see (iii)).

1. (i) $S_x\in \Ore_l(A)$ (since $x^jA_i \subseteq Dx^{i+j}$ for all $j\geq 0$ such that $i+j\geq 0$).

(ii) $\ass_l(S_x) = \ga$: The element $x$ is a homogeneous element of the $\Z$-graded ring $A$. Therefore, $\ga':= \ass_l(S_x) = \oplus_{i\in \Z} \ga_i'$ where $\ga_i':= \ga' \cap A_i$. For $i\geq 0$, $\ga_i' = \CK (\s ) x^i$ since $x^j\cdot : Dx^i\ra Dx^{i+j}$, $ dx^i\mapsto \s^j(d) x^{i+j}$.  Finally, $\ga_{-i}'= \s^{-i} (\CK (\s):(i,-i))y^i$ for $i\geq 1$ since, for all $j\geq 1$,
$$x^{i+j}\cdot : Dy^i\ra Dx^j, \;\; dy^i\mapsto \s^{i+j}(d) (i+j, -i) x^j = \s^{i+j} (d) \s^j((i,-i)) x^i = \s^j (\s^i (d)(i,-i))x^j.$$
Therefore, $\ga' = \ga$.

(iii) $S_x\in \Den_l(R, \ga )$: We have to show that $\gb \subseteq \ga$, or equivalently, $\gb_i\subseteq \ga$ for all $i\geq 1$, by statement 3. If $b\in \gb_{-i}$ then $b(-i,i)=0$. By applying $\s^i$, we get $0=\s^i(b) \s^i ((-i,i))=\s^i (b)(i,-i)$, by (\ref{iist}). Therefore, $b\in \s^{-i}(0:(i,-i))\subseteq \s^{-i}(\CK (\s):(i,-i))\subseteq \ga$.

2. Statement 2 follows from statement 1.

4. The ring $A= \oplus_{i\in \Z}A_i$ is $\Z$-graded and the element $x\in A_1$ is  homogeneous, hence all the elements of the set $S_x$ are homogeneous. Therefore, the ring $A_x= \oplus_{i\in \Z} A_{x,i}$ is automatically the $\Z$-graded ring where $A_{x,0} = \cup_{i\geq 0} x^{-i} D(\s ) x^i$ is the zero component of the ring $A_x$ and $A_{x,i}= A_{x,0}x^i= x^iA_{x,0}$ for all $i\in \Z$.  By (\ref{RsKs}), the map $\overline{\s} : D(\s ) \ra D(\s )$, $\overline{d} := d+\CK (\s )\mapsto \s (d) +\CK (\s )$, is a ring monomorphism. Every element $x^{-i}\overline{d} x^i$ of the ring $A_{x,0}$ where $i\geq 0$ and  $\overline{d} \in D(\s )$, can be written also as follows
%\marginpar{xidxi}
\begin{equation}\label{xidxi}
x^{-i}\overline{d}x^i = x^{-i-j} x^j \overline{d} x^i= x^{-i-j}\overline{\s^j (d)}x^{i+j}\;\; {\rm for}\;\; j\geq 0.
\end{equation}
So, the addition and multiplication in the ring $A_{x,0}$ are given by the rule:
%\marginpar{xidxi1}
\begin{equation}\label{xidxi1}
x^{-i}\overline{d}x^i+ x^{-j}\overline{e}x^j= x^{-i-j}\bigg( \overline{\s^j(d)+\s^i(e)}\bigg) x^{i+j}\;\; {\rm and}\;\;
 x^{-i}\overline{d}x^i\cdot  x^{-j}\overline{e}x^j= x^{-i-j}\overline{\s^j(d)\cdot \s^i(e)}x^{i+j}.
\end{equation}
Since $x\cdot x^{-i}\overline{d}x^i= x^{-i}\overline{\s (d)}x^i\cdot x= \s (x^{-i}\overline{d}x^i)x$, the ring $A_x$ is the skew polynomial ring $A_{x,0}[x^{\pm 1}; \s ]$.

5(a) The equality in the statement (a) follows from the fact that the ring $A$ is a $\Z$-graded ring and the following equalities: For all $d\in D$, $x\cdot dx^i = \s (d) x^{i+1}$ $(i\geq 0)$ and $x\cdot dy^i = \s (da) y^{i-1}$ $(i\geq 1)$. Using the equality in the statement (a), we see that $x\in \CC_A'$ iff $\s$ is a monomorphism and $0=(\ker (\s) :a) = (0:a) = \{ d\in D\, | \, da=0\}$ iff $\s$ is a monomorphism and $a\in \lCCD$.

(b) The equality in the statement (b) follows from the fact that the ring $A$ is a $\Z$-graded ring and the following equalities: For all $d\in D$, $ dx^i\cdot x = d x^{i+1}$ $(i\geq 0)$ and $ dy^i \cdot x= d\tau^{i-1}(a) y^{i-1}$ $(i\geq 1)$.Hence,  $x\in \lCCA$ iff $\tau^i(a)\in \lCCD$ for all $i\geq 0$.

(c) By statement 3, $x\in \CC_A$ iff $x\in \CC_A'$ iff $\s$ is a monomorphism and $a\in \lCCD$.

6. If $\ass_l(S_x) = \ass_r(S_x)$ then, by statements 1 and 3, $\CK (\s ) =0$. Clearly, $\CK (\s ) =0$  iff $\s$ is a monomorphism. Conversely, suppose that $\s$ is a monomorphism. Then $\CK (\s ) =0$. Then, by statements 1 and 3, $\ass_l(S_x) = \ass_r(S_x)$ iff $\s$ is a monomorphism and $\s^{-i}(0:(i,-i)_D)=\ker (\cdot (i,-i)_D)$. An element $d\in D$ belongs to the set $\s^{-i}(0:(i,-i)_D)$ iff $0=\s^i (d) (i,-i)=\s^i (d)\s^i ((-i,i))=\s^i(d(-i,i))$ iff $d(-i,i)=0$ (since $\s$ is a monomorphism) iff $x\in \ker (\cdot (-i,i)_D)$.

7. By statement 1, $S_x\in \Den (A)$ iff $S_x\in \Den_r(S_x)$ (by statement 3) iff $S_x$ is a right Ore set of $A$ and $\s$ is \ monomorphism (by statement 6) iff $\s$ is an automorphism. In more detail, if $\s$ is automorphism then $S_x$ is a right Ore set of $A$ and $\s$ is a monomorphism. Conversely, since $S_x$ is a right Ore set of $A$ then for given elements $d\in D$ and $x\in S_x$, $dx^i=xr$ for some $x^i\in S_x$ and $r\in A$. Without loss of generality  we may assume that $i\geq 1$ (if $i=0$ then $dx=xrx=xr'$ where $r' = rx\in A$). Then $dx^i=x\cdot x^{i-1}d'= \s^i (d') x^i$ for some $d'\in D$ (the ring $A$ is a $\Z$-graded ring). Then $d=\s (\s^{i-1}(d'))$, i.e., $\s$ is an epimorphism. Therefore, $\s$ is an isomorphism (since $\s$ is a monomorphism). $\Box$

Proposition \ref{b1Apr} is a source of  natural examples of `exotic' situations:

(i) {\em A left denominator set which is not a right denominator set}: $S_x$ is a left denominator set of $A$ which is not a right denominator set iff $\s$ is not an automorphism, by statement 7.

(ii) {\em A left denominator set  $S$ such that} $\ass_l(S)\varsupsetneqq \ass_r(S)$: $\ass_l(S)\varsupsetneqq \ass_r(S)$ iff $\s$ is not a monomorphism, by statement 6.

Using (\ref{Asyma1}) or (\ref{Asyma}), the next proposition follows from Proposition \ref{b1Apr}.

\begin{proposition}\label{d2Apr}%\marginpar{d2Apr}
Let $A= D[x,y; \s , \tau , a]$ be a GWA and $S_y=\{ y^i \, | \, i\geq 0\}$. Then
\begin{enumerate}
\item $S_y\in \Den_l(A, \ga^o )$ where $\ga^o := \bigoplus_{i\geq 1} \CK (\tau ) y^i \oplus \bigoplus_{i\geq 0} \tau^{-i} (\CK  (\tau ) :(-i,i)) x^i$ and $D(\tau ) := D/\CK (\tau )$ (recall that $o$ is the symmetry isomorphism of $A$).
\item $A/ \ga^o \simeq \bigoplus_{i\geq 1} D (\tau )y^i \oplus \bigoplus_{i\geq 0}  D/\tau^{-i} (\CK  (\tau ) :(-i,i))  x^i$  and the map $A/\ga \ra A/ \ga^o$, $u+\ga \mapsto u^o+\ga^o$, is a ring isomorphism.
\item $\ass_l(S_y)=\bigoplus_{i\geq 1} \ker (\cdot (i,-i)_D)x^i\subseteq \ga^o$.
\item $A_y:= S_y^{-1}A=A_{y,0}[y^{\pm 1}; \tau ]$ is a skew polynomial ring where $A_{y,0}:= \bigcup_{i\geq 1} y^{-i}D (\tau )y^i$ and $\tau (y^{-i} (d+\CK (\tau )) y^i) := y^{-1} (\tau (d)+\CK (\tau )) y^i$. The addition and multiplication in the ring $A_{y,0}$ is given in a similar fashion as in (\ref{xidxi1}).
   % \item The map $o : A_x\mapsto A_y$, $x^{-i} b \mapsto (x^{-i}b)^o= y^{-i} b^o$, is an isomorphism of $\Z$-graded rings, i.e.,  $A_{x,i}^o=A_{y,i}$ for all $i\in \Z$.
   \item
  \begin{enumerate}
\item $\ker_A(y\cdot ) =\bigoplus_{i\geq 0}\ker (\tau ) y^i\oplus
\bigoplus_{i\geq 1}(\ker (\tau ) : \s (a))x^i$; $y\in \CC_A'$ iff $\tau$ a monomorphism  and $\s (a)\in \lCCD $.
\item $\ker_A(\cdot y) = \bigoplus_{i\geq 1}(\ker_D (\cdot \s^i(a))x^i$; $y\in \lCCA$ iff $\s^i(a)\in \lCCD$ for all $i\geq 1$.
\item $y\in \CC_A$ iff $y\in \CC_A'$ iff $\tau$ is a monomorphism and $\s (a)\in \lCCD$.
\end{enumerate}
\item $\ass_l(S_y)=\ass_r(S_y)$ iff $\tau$ is a monomorphism.
\item $S_y\in \Den (R)$ iff $\tau$ is an automorphism. If $S_y\in \Den (R)$ then $\ass (S_y) =\bigoplus_{i\geq 1} \ker (\cdot (i,-i)_D)x^i$.
\end{enumerate}
\end{proposition}

{\bf A regularity criterion for the elements $x$ and $y$ of a GWA.} Such a criterion is the following proposition which is used in the proof of a simplicity criterion for GWAs (Theorem \ref{29Oct}).

\begin{proposition}\label{a30Oct}%\marginpar{a30Oct}
The following statements are equivalent.
\begin{enumerate}
\item $x,y\in \CC_A$.
\item $\s$ and $\tau$ are monomorphisms  and $a, \s (a) \in \lCCA$.
\item $a, \s (a) \in \CC_D$.
\end{enumerate}
\end{proposition}

{\it Proof}. $(1\Leftrightarrow 2)$ The equivalence follows from Proposition \ref{b1Apr}.(5c) and Proposition \ref{d2Apr}.(5c).

$(1\Rightarrow 3)$ If $x,y\in \CC_A$ then $a=yx\in \CC_A$ and $\s (a) = xy\in \CC_A$. Hence, $a, \s (a) \in \CC_D$.

$(3\Rightarrow 1)$ If $a, \s (a) \in \CC_D$ then $a, \s (a) \in \lCCD$ and $\s$, $\tau$ are monomorphisms since $ad=\tau\s (d) a$ and $\s (a)d=\s \tau (d) \s (a) $ for all $d\in D$. $\Box $

%$\noindent $

{\bf A criterion for a GWA to be a domain.}
\begin{proposition}\label{a8Apr}%\marginpar{a8Apr}
Let $A=D(\s, \tau , a)$ be a GWA. Then $A$ is a domain iff $D$ is a domain, $\s$ and $\tau $ are monomorphisms of $D$ and $a\neq 0$.
\end{proposition}

{\it Proof}. $(\Rightarrow)$ If $A$ is a domain then $D$ is a domain and $a\neq 0$ (since $a=yx$). By Proposition \ref{b1Apr}.(5) and Proposition \ref{d2Apr}.(6),
 $\s$ and $\tau $ are monomorphisms of $D$.

$(\Leftarrow)$ Since  $a\neq 0$, all the elements $(n,m)$ (where $n,m\in \Z$) are nonzero. Since $D$ is a domain and $\s$ and $\tau$ are monomorphisms of $D$, the ring $A$ is a domain, by (\ref{dmult1}).  $\Box $

{\bf The ideals $(x^n)$ and $(y^n)$ of GWAs}.
For a ring $R$ and an element $r\in R$, $(r) := RrR$ is the ideal of $R$ generated by the element $r$. The next lemma describes the ideals $(v_n)$ of a GWA $A=\oplus_{i\in \Z} Dv_i$ where $n\in \Z$. Clearly, $(v_0)=(1)=A$.

\begin{lemma}\label{a5Apr}%\marginpar{a5Apr}
Let $A=D(\s, \tau , a)$ be a GWA. Then, $(v_n) = \oplus_{m\in \Z} (v_n)_m$ where $(v_n)_m:= (v_n)\cap Dv_m$ %$= Dv_nv_{m-n} +Dv_{m-n} v_n = (D(n,m-n) D+(m-n, n))v_m$.
For all $n\geq 1$,
$$(x^n)_m=\begin{cases}
Dx^m& \text{if }m\geq n,\\
\sum_{i=0}^{n-m}Dy^ix^ny^{n-m-i}=\sum_{i=0}^{n-m}D(-i,n)(n-i, m-n+i)v_m& \text{if }m<n, \\
\end{cases}$$

$$(y^n)_m=\begin{cases}
Dy^{-m}& \text{if }m\leq -n,\\
\sum_{i=0}^{n+m}Dx^iy^nx^{n+m-i}=\sum_{i=0}^{n+m}D(i,-n)(-n+i, m+n-i)v_m& \text{if }m>-n. \\
\end{cases}$$
\end{lemma}

{\it Proof}. The GWA $A= \oplus_{i\in \Z} A_i$ is a $\Z$-graded ring and the element $v_n$ is a homogeneous element of $A$. Hence, $(v_n) = \oplus_{m\in \Z} (v_n)_m$ is a homogeneous ideal of $A$, i.e.,  $(v_n)_m = (v_n)\cap A_m$. Now, $(v_n)_m = \sum_{i+j+n=m} A_iv_nA_j=\sum_{i+j+n=m}Dv_iv_nDv_j=\sum_{i+j+n=m} Dv_iv_nv_j$. Let $n\geq 1$, that is $v_n=x^n$. Then $(v_n)_m= Dx^m$ for all $m\geq n$. Let $m<n$. Then $\D = n-m>0$. The product $v_ix^nv_{-\D -i}$ is equal to
\begin{eqnarray*}
i>0: & & x^ix^ny^{\D +i} = x^n (x^iy^i) y^\D = \s^n((i,-i))x^ny^\D,  \\
-\D \leq i \leq 0: & & y^{-i}x^ny^{\D +i}=(i,n) (n+i, m-n-i)v_m,  \\
i<-\D : & & y^{-i}x^nx^{-\D -i}=y^\D (y^{-\D -i}x^{-\D -i} ) x^n = \tau^{-\D}((\D +i, -\D -i))y^\D x^n.  \\
\end{eqnarray*}
Therefore, $(x^n)_m=\sum_{-\D \leq j \leq 0} Dy^{-j}x^n y^{-\D -j}  =\sum_{i=0}^{n-m} Dy^ix^ny^{n-m-i} = \sum_{i=0}^{n-m}D(-i,n) (n-i, m-n+i) v_m$. By the canonical $(x,y)$-symmetry isomorphism, for $n\geq 1$, $(y^n)_m=Dy^{-m}$ for all $m\leq -n$, and, for $m>-n$, $(y^n)_m= \sum_{i=0}^{m+n}Dx^iy^nx^{m+n-i}= \sum_{i=0}^{m+n} D(i,-n) ( -n+i, m+n-i)v_m$. $\Box$

%Using (\ref{1vnvmm}) and explicit expressions for $(i,j)$, we see that $(v_n)_m= (D(n,m-n)+D(m-n,n))v_m = Dv_nv_{m-n} +Dv_{m-n} v_n$.

The next corollary is used in the proof of a simplicity criterion for GWAs (Theorem \ref{B5Apr}).
\begin{corollary}\label{c5Apr}%\marginpar{c5Apr}
Let $A= D[x,y; \s , \tau , a]$ be a GWA and $n\geq 1$  a natural number. Then
\begin{enumerate}
\item $(x^i)=A$ for $i=1, \ldots , n$ iff $Da+D\s^i(a) =D$ for $i=1, \ldots , n$. In particular,  $(x^i) = A$ for all $i\geq 1$  iff $Da+D\s^i(a) = D$ for all $i\geq 1$.
\item $(y^i)=A$ for $i=1, \ldots , n$ iff $D\s (a)+D\tau^i\s (a) =D$ for $i=1, \ldots , n$. In particular, $(y^i) = A$ for all $i\geq 1$  iff $D\s (a)+D\tau^i\s (a) = D$ for all $i\geq 1$.
\end{enumerate}
\end{corollary}

{\it Proof}. 1. $(\Rightarrow)$ Let $\ga_i := Da+D\s^i (a)$. By Lemma \ref{a5Apr}, for all $i\geq 1$, $(x^i)_{i-1} = (D(i,-1) +D(-1, i))x^{i-1} = (D\s^i(a) + Da)x^{i-1} = \ga_i x^{i-1}$. If $(x^i) = A$ then $\ga_i = D$.

$(\Leftarrow )$ Suppose that $\ga_i =D$ for all $i=1, \ldots , n$. We use induction on $n$ to show that $(x^i) = A$ for all $i=1, \ldots , n$. For $n=1$, $(x) = A$ since $(x)_0= \ga_1=D$. Suppose that $n>1$ and $(x) = \cdots =(x^n)=A$. Since $(x^{n+1})_n= \ga_{n+1} x^n = Dx^n$, we have the inclusion $(x^{n+1})\supseteq (x^n) = A$, i.e., $(x^{n+1}) = A$.

2. Statement 2 follows from statement 1 by the $(x,y)$-symmetry. $\Box $

{\it Example.} Let $A=D[x,y; \s, \s^{-1}, a=h(h-n)]$ where $D= K[h]$, $\s (h) = h-1$ and $n\geq 1$ is a natural number. By Corollary \ref{c5Apr}.(1), $(x) = (x^2)=\cdots = (x^{n-1})=A$ and $(y) = (y^2)=\cdots = (y^{n-1})=A$ (since $Da+D\s^i(a) = D$ and $D\s (a)+D\s^{-i}\s(a)=\s^{-i+1}(Da+D\s^i(a)) = D$ for $i=1, \ldots , n-1$) and $(x^n) \neq A$ and $(y^n) \neq A$ since $Da+D\s^n(a) = D\s^n(h) \neq D$ and $D\s (a) +D\s^{-n}\s(a) = \s^{-n+1}(D\s^n(h))= D\s (h) \neq D$.

In case when the elements $a$ and $\s (a)$ are normal and regular, the conditions (\ref{abGWA}) can be simplified, as the following lemma shows. It also gives a criterion for the endomorphisms $\s$ and $\tau$ to commute.

\begin{lemma}\label{b5Apr}%\marginpar{b5Apr}
Let $A= D(\s , \tau , a)$ be a GWA.  Suppose that $a$ and $\s (a)$ are  normal,  regular elements in $D$.  Then
\begin{enumerate}
\item The conditions (\ref{abGWA}) are equivalent to the equalities $\tau \s = \o_a$ and $ \s \tau = \o_{\s (a)}$.
\item $\s$ and $\tau$ are automorphisms of the ring $D$ and $\tau = \o_a\s^{-1}$.
\item $\s \tau = \tau \s$ iff $ \o_a= \o_{\s (a)}$ iff there is a normal regular element $b\in D$ such that $ab= ba$, $a^i \s (a) = ab$ and $\o_b = \o_{a^i}$ for some $i\geq 0$.
\end{enumerate}
\end{lemma}

{\it Proof}. 1. Statement 1 is obvious.

2. Statement 2 follows from statement 1 since $\o_a$ and $\o_{\s (a)}$ are automorphisms, by Lemma \ref{a2Apr} (as the elements $a$ and $\s (a)$ are normal and regular).

3. By statement 1, $\s \tau = \tau \s$ iff $ \o_a= \o_{\s (a)}$.

 Suppose that  $ \o_a= \o_{\s (a)}$. The element $a$ of $D$ is normal and regular. Then the ring $D$ is a subring of the localization $D_a$ of $D$ at the powers of the element $a$. Since $\o_a(a) =a$, the automorphism $\o_a$  of $D$ is extended in a unique way to an automorphism of $D_a$: $\o_a(a^{-i}d) =a^{-i} \o_a(d)$ where $i\geq 0$ and $d\in D$. Then the equality $\o_a = \o_{\s (a)}$ implies  $\o_{a^{-1} \s (a)} = \id_{D_a}$, and so  $a^{-1} \s (a) = z$ for some central regular element $z\in D_a$. So, $z= a^{-i} b$ for some $i\geq 0$ and $b\in D$ or, equivalently, $b= a^i z$. Then $ab=ba$.  The element $b\in D$ is normal and regular and $\o_b = \o_{a^i}$. The equality $z= a^{-i} b$ can be written as $a^i \s (a) = ab$ (by using $z= a^{-1} \s (a)$).

Conversely, if $ ab=ba$, $a^i \s (a) = ab$ and $\o_b = \o_{a^i}$ for some normal regular element $b\in D$ and $i\geq 0$. Then the equality $a^i \s (a) = ba$ implies the equality $\o_{a^i} \o_{\s (a)} = \o_b \o_a= \o_{a^i} \o_a$, hence $ \o_{\s (a)} = \o_a$. $\Box$

%%%%%%%%%%%%%%%%%% SECTION 3  %%%%%%%%%%%%%%%%%%%%%%%%

\section{Simplicity criteria for Generalized Weyl algebras}\label{SIMCGWA}%\marginpar{SIMCGWA}

The aim of this section is to give two simplicity criteria for GWAs (Theorem \ref{B5Apr} and Theorem \ref{29Oct}). The first one  (Theorem \ref{B5Apr}) is a simplicity criterion for a GWA $A=D[x,y; \s , \tau , a]$ where the elements $a$ and $\s (a)$ are normal in $D$. This is a mild restriction on the elements $a$ and $\s (a)$ that often occurs in applications.  The second one (Theorem \ref{29Oct}) is
 a simplicity criterion for GWAs in general case (no restrictions on $a$ and $\s (a)$). Their proofs are quite different.

{\bf Simplicity criteria via denominator sets.}
Let $R$ be a ring, $S\in \Den_l(R)$ and $T\in \Den_r(R)$. Let $I$ be an ideal of $R$. In general, neither the left ideal $S^{-1}I$ of $S^{-1}R$ nor the right ideal $IT^{-1}$ of $RT^{-1}$ is an ideal. But if the ring $R$ is commutative or Noetherian then one-sided ideals $S^{-1}I$ and  $IT^{-1}$ are ideals of the localized rings. So, in general, ideals of a ring and its localization are not much related.   Suppose that the rings $S^{-1}R$ and $RT^{-1}$ are $R$-{\em isomorphic}, i.e., there is a ring isomorphism $f: S^{-1}R\ra RT^{-1}$ such that $f(r \alpha ) = rf(\alpha )$ for all elements $r\in R$ and $\alpha \in S^{-1}R$. In particular, $\ass_l(S) = \ass_r(T)$. We identify the rings $S^{-1}R$ and $RT^{-1}$.

{\it Example}. Let $S$ be a left and right denominator set of  a ring $R$. Then the rings $S^{-1}R$ and $RS^{-1}$ are $R$-isomorphic.

The third statement of the following proposition is a useful simplicity criterion for a ring via its localization at a {\em left} and {\em right} denominator set. This criterion is used in the proof of a simplicity criterion (Theorem \ref{B5Apr}). The fourth statement is a simplicity criterion for a ring via its localizations at a left denominator set and a right denominator set.

\begin{proposition}\label{A5Apr}%\marginpar{A5Apr}
Let $R$ be a ring, $S\in \Den_l(R)$ and $T\in \Den_r(R)$ be such that the rings $S^{-1}R$ and $RT^{-1}$ are $R$-isomorphic and $I:=\ass_l(S)= \ass_r(T)$. Let $\ga$ be an ideal of $R$ and $S^{-1}R\ga RT^{-1}= S^{-1} \ga T^{-1}$ be the ideal of $S^{-1}R=RT^{-1}$ generated by $\ga$.  Then
\begin{enumerate}
\item $S^{-1} \ga T^{-1}=S^{-1}R$ iff $\ga \cap ST\neq \emptyset$ where $ST:= \{ st\, | \, s\in S, t\in T\}$.
\item If $S\in \Den (R)$ then $S^{-1} \ga S^{-1}=S^{-1}R$ iff $\ga \cap S\neq \emptyset$.
\item If $S\in \Den (R)$ then $R$ is a simple ring iff $S^{-1}R$ is a simple ring, $\ass (S)=0$ and $RsR=R$ for all $s\in S$.
    \item The ring $R$ is simple iff $S^{-1}R=RT^{-1}$ is a simple ring, $I=0$ and $RstR=R$ for all elements $s\in S$ and $t\in T$.
\end{enumerate}
\end{proposition}

{\it Proof}. 1. $S^{-1}\ga T^{-1} = S^{-1}R$ iff $1=\sum_{i=1}^n s_i^{-1} a_it_i^{-1}$ for some elements $s_i\in S$ and $t_i\in T$. Choose elements $s\in S$ and $t\in T$ such that  $ss_i^{-1} = \frac{r_i}{1}$ and $t_i^{-1} t = \frac{r_i'}{1}$ for some elements $r_i, r_i'\in R$. Then $st = \sum_{i=1}^n \frac{r_i}{1}a_i \frac{r_i'}{1}= \frac{a}{1}$ where $a=\sum_{i=1}^n r_ia_i r_i\in \ga$, Hence, $s'st = s'a\in ST\cap \ga$ for some element $s'\in S$.

2. Statement 2 follows from statement 1.

3. Statement 3 follows from statement 2.

4. Statement 4 follows from statement 1.  $\Box $

{\bf Simplicity criteria for generalized Weyl algebras}.{\bf Proof of Theorem \ref{B5Apr}}. $(1\Rightarrow 2)$ Suppose that $A$ is a simple ring. By Proposition \ref{b1Apr} and Proposition \ref{d2Apr}, $\ass_l(S_x)=0$ and $\ass_l(S_y)=0$, i.e., the elements $x$ and $y$ are regular in $A$. Then the elements  $a=yx$ and $\s (a) = xy$ are regular in $A$. In particular, they are regular in $D$, and so that the statement (a)  holds.

The elements $a$ and $\s (a)$ are normal and regular. By Lemma \ref{b5Apr}, $\s$ and $\tau$ are automorphisms of $D$. Hence, $S_x$ (and $S_y$) is a left and right Ore set in $A$. By Proposition \ref{A5Apr}.(3), the ring  $A$ is  simple iff the ring  $S_x^{-1} A$ is simple and $(x^i) = A$ for all $i\geq 1$. It is a classical result that the ring $S_x^{-1}A= D[x^{\pm 1}; \s ]$ is simple iff the conditions (b) and (c) hold. By Corollary \ref{c5Apr}.(1), $(x^i)=A$ for all $i\geq 1$ iff the condition (d) holds.

$(2\Rightarrow 1)$ Suppose that the conditions (a)-(d) of statement 2 hold. By the assumption, the elements $a$ and $\s (a)$ are normal, and they are regular, by the statement (a). By Lemma \ref{b5Apr}, $\s$ and $\tau$ are automorphisms of the ring $D$. Therefore, the elements $(n,m)$ (where $n,m\in \Z$) are regular, see (\ref{1vnvmm}). Then the elements $x$ and $y$ are regular. Hence, $S_x$ (and $S_y$) is a left and right Ore set in $A$. Then repeating the arguments at the end of the proof of the implication $(1\Rightarrow 2)$, which are of `iff'-nature,   we see that $A$ is a simple ring (By Proposition \ref{A5Apr}.(3), the ring  $A$ is  simple iff the ring  $S_x^{-1} A$ is simple and $(x^i) = A$ for all $i\geq 1$. It is a classical result that the ring $S_x^{-1}A= D[x^{\pm 1}; \s ]$ is simple iff the conditions (b) and (c) hold. By Corollary \ref{c5Apr}.(1), $(x^i)=A$ for all $i\geq 1$ iff the condition (d) holds.).

$(1\Leftrightarrow 3)$  In view of the $(x,y)$-symmetry isomorphism of the GWA $A$, the equivalence $(1\Leftrightarrow 3)$ follows from the equivalence $(1\Leftrightarrow 2)$.  $\Box $

{\em Example.} Let $A= \L [x,y; \s = \s_{\alpha , \beta}, \tau = \o_a\s^{-1}, a=zp^iq^j]$ be the GWA considered in Section \ref{GWA2}. If $i+j>0$ then $a,\s^i (a) \in (p^iq^j)\neq D$ for all $i\geq 1$, and so the condition (2d) of Theorem \ref{B5Apr} does not hold and the algebra $A$ is not simple.

For a ring $D$ and its ring endomorphism $\s$, the subring of $D$, $D^\s=\{ d\in D\, | \, \s (d) = d\}$,   is called the {\em ring of} $\s$-{\em invariants}, and each element of $D^\s$ is called a $\s$-{\em invariant}.  Every left normal, left regular element $d$ of $D$, determines a ring endomorphism of $D$:
%\marginpar{lddp}
\begin{equation}\label{lddp}
\o_d:D\ra D, \;\; d'\ra \o_d(d'), \;\; {\rm where}\;\; dd' = \o_d(d')d.
\end{equation}
The next theorem is a simplicity criterion for GWAs.

\begin{theorem}\label{29Oct}%\marginpar{29Oct}
 Let $A= D[x,y; \s, \tau , a]$ be a GWA. Then the following statements are equivalent.
\begin{enumerate}
\item  $A$ is a simple ring.
\item
\begin{enumerate}
\item The elements $a$ and $\s (a)$ are regular in $D$,
\item For all nonzero ideals $I$ of $D$, $I'=D$ where
$I':= I+\sum_{i\geq 1} \Big( D\s^i(I)(i,-i)+D\tau^i(I)(-i,i)\Big)$.
\item None of the ring endomorphisms   $\s^n$ $(n\geq 1)$ of $D$ is equal to the ring endomorphism $\o_d$ (see (\ref{lddp})) where $d$ is a $\s$-invariant, regular, left normal element of $D$.
\end{enumerate}
\item
\begin{enumerate}
\item The elements $a$ and $\s (a)$ are regular in $D$,
\item For all nonzero ideals $I$ of $D$, $I'=D$ where
$I':= I+\sum_{i\geq 1} \Big(D\s^i(I)(i,-i)+D\tau^i(I)(-i,i)\Big)$.
\item None of the ring endomorphisms   $\tau^n$ $(n\geq 1)$ of $D$ is equal to the ring endomorphism $\o_d$ (see (\ref{lddp})) where $d$ is a $\tau$-invariant, regular, left normal element of $D$.
\end{enumerate}
\end{enumerate}
If one of the equivalent conditions holds then $\s$ and $\tau$ are monomorphisms of $D$.
\end{theorem}

{\it Proof}. $(1\Rightarrow 2)$ Suppose that $A$ is a simple ring. By Proposition \ref{b1Apr} and Proposition \ref{d2Apr}, $\ass_l(S_x)=0$ and $\ass_l(S_y)=0$, i.e., the elements $x$ and $y$ are regular in $A$. Then the elements  $a=yx$ and $\s (a) = xy$ are regular in $A$. In particular, they are regular in $D$, and so that the statement (a)  holds.

Let $I$ be a nonzero ideal of $D$. The zero component $(I)_0$ of the ideal $$(I)= AIA= \sum_{i,j\in \Z} Dv_iIDv_j= \sum_{i,j\in \Z} D\s^i(I) v_iv_j=\sum_{i,j\in \Z} D\s^i(I) (i,j) v_{i+j}$$
 is equal to $I'=\sum_{i\in \Z} D\s^i(I)(i,-i)$. This is the LHS of the equality in the statement (b). The algebra $A$ is simple and the ideal $(I)$ is a nonzero homogeneous ideal of $A$. Therefore, $(I)_0=D$, i.e., the statement (b) holds.

 Finally, suppose that the statement (c) is false, i.e., $\s^n = \o_d$ for some $\s$-invariant, regular, left normal element of $D$. Then $dd'= \s^n (d') d$ for all elements $d'\in D$. We claim that {\em the ideal of $A$ generated by the element $u=x^n+d$ is not equal to} $A$. Notice that $ux=xu$  (since $\s (d) = d$) and $ud' = \s^n(d')u$ for all $d'\in D$ (since $\s^n = \o_d$). This means that the element $u\in A_+:=\oplus_{i\geq 0} Dx^i$ is a left normal element of the ring $A_+$. In particular, $A_+uA_+= A_+u$ is an ideal of $A_+$. Suppose that $AuA=A$, i.e., $1\in \sum_{-l\leq i,j\leq l} A_iuA_j$ for some $l\in \N$, we seek a contradiction. Then $x^{2l}= x^l1x^l \in A_+uA_+= A_+u$, and so $x^{2l} = vu$ for some element $v\in A_+$. Then $v= x^{2l-n}+\cdots + d_mx^m$ where $d_m\in D\backslash \{ 0\}$ and $d_mx^m$ is the least term of $v$ (w.r.t. the $\Z$-grading of $A$). Then $0=d_mx^md= d_m\s^m(d) x^m = d_mdx^m$ (since $\s (d) =d$), i.e., $d_md=0$. This contradicts to regularity of $d$. This means that $0\neq (u) \neq A$, as claimed. Therefore, the statement (c) holds.

$(2\Rightarrow 1)$ Suppose that the conditions (a) - (c) hold. By the statement (a), the elements $a$ and $\s (a)$ are regular in $D$. By Proposition  \ref{a30Oct},  $x, y\in \CC_A$ and $\s , \tau$ are monomorphisms. Hence, $x^i, y^i \in \CC_A$ for all $i\geq 1$.

Let $J$ be a nonzero ideal of $A$. We have to show that $J=A$. The ideal $J$ contains a nonzero element, say $u$, of least possible length, say $l$. If $l=0$ then $J$ contains a nonzero element, say $d$, from $D$ (since $x^i$ and $y^i$ are regular element of $A$ for all $i\geq 1$). Let $I= DdD$. Then $J\supseteq I+\sum_{i\geq 1} (Dx^iIy^i+Dy^iIx^i)= I' =D$, by the statement (b), and so $J=A$, sa required.

Suppose that $l>1$. Replacing the element $u$ by $x^su$ or $y^su$ for some $s\geq 0$, we may assume that
 $u=u_0+u_1y+\cdots +u_ly^l$ for some elements  $u_i\in D$ such that  $u_0\neq 0$ and $u_l\neq 0$ (since $x^i$ and $y^i$ are regular elements in $A$). Let $I=Du_0D$, a nonzero ideal of $D$. Then
 $$ J\supseteq DuD+\sum_{i\geq 1} (Dx^iDuDy^i+Dy^iDuDx^i) = I+\sum_{i\geq 1} (D\s^i(I) (i,-i) + D\tau^i(I) (-i,i))+\cdots  = I'+\cdots  $$
where the three dots means smaller terms, i.e.,  elements of the set $\oplus_{i\geq 1} Dy^i$. By the statement (b), the ideal $J$ contains an element of the form $v=1+\cdots$. Then $0\neq  w:= vx^l= d_0+d_1x+\cdots + d_{l-1}x^{l-1}+x^l\in J$ where $d_i\in D$ and $d_0\neq 0$. By the minimality of $l$, the element $[x,w]= \sum_{i=0}^{l-1}(\s (d_i) - d_i) x^{i+1} \in J$ must be zero, i.e., $xw=wx$ and $\s (d_i) = d_i$ for all $i=0, 1, \ldots , l-1$. Similarly, $\s^l(d) w-wd= \sum_{i=0}^{l-1}(\s^l (d) d_i-d_i\s^i (d))x^i \in J$ for all $d\in D$. Therefore, $\s^l (d) w=wd$ and $\s^l (d) d_i = d_i\s^i(d)$ for all $i=0, 1, \ldots , l-1$. In particular, $\s^l (d) d_0= d_0d$ for all $d\in D$. The element $d_0$ is a regular element of $D$. Since otherwise we would have either $d'd_0=0$ or $ d_0d'=0$ for some $d'\in D$. Then either $0\neq d' w= \sum_{i=1}^{l-1} d' d_ix^i +d' x^l \in J$ or $0\neq wd' = \sum_{i=1}^{l-1} d_i\s^i(d') x^i  +\s^l(d') x^l \in J$. In both cases, this  would contradict  the minimality of $l$. Therefore, $J=A$, as required.

$(1\Leftrightarrow 3)$ This equivalence follows from the equivalence
 $(1\Leftrightarrow 2)$ by the $(x,y)$-symmetry.  $\Box$

{\bf Involutions on GWAs}. An anti-isomorphism $*$ of a ring $R$ ($(ab)^*= b^*a^*$ for all $a,b\in R$) is called an {\em involution} if $a^{**}=a$ for all elements $a\in R$.

The Weyl algebra $A_1$ admits the canonical involution $*$: $x^*=\der$ and $\der^*=x$. Recall that the Weyl algebra $A_1$ is the GWA $K[h][x, \der ; \s, \s^{-1}, a=h]$ and the involution $*$ respects the subalgebra $K[h]$: $K[h]^*=K[h]$ since $h^*=h$. So, the involution $*$ on $A_1$ can be seen as an extension of the (trivial) involution on the commutative algebra $K[h]$ to $A_1$. The following lemma and its corollary explore further this fact/idea.

\begin{lemma}\label{b2Apr}%\marginpar{b2Apr}
Let $A=D[x,y; \s , \tau , a]$ be a GWA. Suppose that $*$ is an involution of the ring $D$  such that $\s *\tau = *$, $\tau *\s = *$, $a^*=a$ and $\s (a) ^*= \s (a)$. Then the involution $*$ can be extended to an involution $*$ of $A$ by the rule $x^*=y$ and $y^*=x$.
\end{lemma}

{\it Proof}. Let us show that $*$ respects the conditions (\ref{abGWA}) of the GWA $A$: Let $d\in D$.

(a) $ad= \tau \s (d) a$: $ (\tau \s (d) a)^* = a^* \cdot *\tau \s (d) = a\cdot *\tau \s (d) = \tau (\s *\tau ) \s (d) a= \tau *\s (d) a = d^* a^* = (ad)^*$.

(b) $\s (a) d= \s \tau (d) \s (a)$: $ (\s \tau (d) \s (a))^*=\s (a)^* \cdot *\s\tau (d) = \s (a) \cdot *\s \tau (d) = \s (\tau *\s ) \tau (d) \s (a) = \s * \tau (d) \s (a) = d^* \s (a)^* = (\s (a)d)^*$.

 Let us show that   $*$ respects the defining relations (\ref{GWADEF}) of the GWA $A$: Let $d\in D$.

(i) $xd=\s (d) x$: $(\s (d) x)^*= y\s (d)^* = \tau *\s (d)y = d^*y= (xd)^*$.

(ii) $yd= \tau (d) y$: $(\tau (d) y)^* = x\tau (d)^* = \s *\tau (d) x = d^*x= (yd)^*$.

(iii) $yx=a$: $(yx)^*= yx = a= a^*$.

(iv) $xy= \s (a)$: $(xy)^*= xy= \s (a) = \s (a)^*$. $\Box $

For a {\em commutative} ring $D$, the identity map of $D$ is an involution on $D$ which is called the {\em trivial involution} on $D$.

\begin{corollary}\label{c2Apr}%\marginpar{c2Apr}
Let $A=D[x,y; \s , \tau , a]$ be a GWA where $D$ is a commutative ring, $\s$ and $\tau$ are automorphisms of $D$ such that $\tau = \s^{-1}$. Then the trivial involution on $D$ can be extended to an involution $*$ of $D$ by the rule  $x^*=y$ and $y^*=x$.
\end{corollary}

{\it Proof}. The result follows from Lemma \ref{b2Apr}. $\Box$

{\it Example.} The universal enveloping algebra of the Lie algebra $\sl_2$  over a field $K$ of characteristic 0, $U= U(\sl_2) = K\langle X,Y,H\, | \, [H,X]=X,\ [H,Y]=-Y,\ [X,Y]=2H\rangle $,  is isomorphic to the classical GWA $K[H,C][X,Y; \sigma,a=C-H(H+1) ],$
where $\sigma :H\rightarrow H-1, C\rightarrow C$, and $C= YX-H(H+1)$ is the {\em Casimir} element of the algebra $U$. As a universal enveloping algebra, the algebra $U$ admits the {\em canonical involution} given by the rule $X\mapsto -X$, $Y\mapsto -Y$, $H\mapsto -H$. Notice that $C\mapsto C$.

 By Corollary \ref{c2Apr}, the  identity involution on the commutative algebra $K[H, C]$ can be extended to the involution $*$ on $U$ given by the rule $X^*=Y$, $Y^*=X$ and $H^* = H$. This is not the canonical involution on $U$.

%%%%%%%%%%%%%%%%%% SECTION 4  %%%%%%%%%%%%%%%%%%%%%%%%

\section{Generalized Weyl algebras of rank $n$}\label{GWARNN}%\marginpar{GWARNN}

The aim of this section is to introduce a new class of rings which is more general that the class of generalized Weyl algebras. The rings of the new class are also called generalized Weyl algebras. In order to distinguish these new rings from the old ones the latter are called {\em classical} GWAs.

{\bf Iterated generalized Weyl algebras.} The next corollary follows from Theorem \ref{1Apr} by induction on the rank $n$.

\begin{corollary}\label{ax1Apr}%\marginpar{ax1Apr}
Let $A= D[x_1,y_1; \s_1, \tau_1, a_1] \ldots [x_n,y_n; \s_n, \tau_n, a_n]$ be an iterated GWA of rank $n$. Then $A= \oplus_{\alpha \in \Z^n} Dv_\alpha$ is a direct sum of the free left $D$-modules ${}_DDv_\alpha \simeq D$ where for $\alpha = (\alpha_1, \ldots , \alpha_n)$, $v_\alpha = v_{\alpha_1}(1)\cdots v_{\alpha_n}(n)$ and $v_{\alpha_i}(i) = \begin{cases}
x_i^{\alpha_i}& \text{if }\alpha_i\geq 0,\\
y_i^{-\alpha_i}& \text{if }\alpha_i<0.\\
\end{cases} $
\end{corollary}

{\bf Classical generalized Weyl algebras, \cite{Bav-GWA-FA-91}-\cite{Bav-Bav-Ideals-II-93}.} Let $D$ be a ring, $\sigma=(\sigma_1,...,\sigma_n)$   an $n$-tuple of
commuting automorphisms of $D$,  $a=(a_1,...,a_n)$   an $n$-tuple  of elements of  the centre
$Z(D)$  of $D$ such that $\sigma_i(a_j)=a_j$ for all $i\neq j$. The {\bf (classical) generalized Weyl algebra} $A=D(\sigma,a)=D[x, y; \sigma,a]$ of rank  $n$  is  a  ring  generated
by $D$  and    $2n$ indeterminates $x_1,...,x_n, y_1,...,y_n$
subject to the defining relations:
$$y_ix_i=a_i,\;\; x_iy_i=\sigma_i(a_i),\;\; x_id=\sigma_i(d)x_i,\;\;{\rm and}\;\;  y_id=\sigma_i^{-1}(d)y_i\;\; {\rm for \; all}\;\; d \in D,$$
$$[x_i,x_j]=[x_i,y_j]=[y_i,y_j]=0, \;\; {\rm for \; all}\;\; i\neq j,$$
where $[x, y]=xy-yx$. We say that  $a$  and $\sigma $ are the  sets  of
{\it defining } elements and automorphisms of the GWA $A$, respectively.

{\em Example.}  The $n$'th {\em Weyl algebra},
 $A_n=A_n(K)$ over a field (a ring) $K$ is  an associative
 $K$-algebra generated by  $2n$ elements
 $x_1, ..., x_n,y_1,..., y_n$, subject to the defining relations:
$$[y_i,x_i]=\d_{ij}\;\; {\rm and}\;\;  [x_i,x_j]=[y_i,y_j]=0\;\; {\rm for\,all}\; i,j, $$
where $\d_{ij}$ is the Kronecker delta function.
The Weyl algebra $A_n$ is a generalized Weyl algebra
 $A=D[x, y; \s ;a]$ of rank $n$ where
$D=K[H_1,...,H_n]$ is a polynomial ring   in $n$ variables with
 coefficients in $K$, $\s = (\s_1, \ldots , \s_n)$ where $\s_i(H_j)=H_j-\delta_{ij}$ and
 $a=(H_1, \ldots , H_n)$.  The map
$$A_n\ra A,\;\; x_i\mapsto  x_i,\;\; y_i \mapsto  y_i,\;\;  i=1,\ldots ,n,$$
is an algebra  isomorphism (notice that $y_ix_i\mapsto H_i$). This was the reason why I called the algebras A above GWAs.

{\bf Generalized Weyl algebras.} Let $A$ be a ring and $\s$ its endomorphism. A subring $B$ of $A$ is called $\s$-{\em invariant} if $\s (B)\subseteq B$.

{\em Definition.} An  iterated generalized Weyl algebra $A=D[x_1,y_1; \s_1 , \tau_1 , a_1]\ldots [x_n, y_n; \s_n , \tau_n , a_n]$ is called a {\bf generalized Weyl algebra} of rank $n$ if $a_1, \ldots , a_n\in D$, the ring $D$ is $\s_i$- and $\tau_i$-invariant for all $i=1, \ldots , n$; and
 for all integers $i,j=1, \ldots , n$ such that $i>j$:
$$\s_i(x_j)=\l_{ij}x_j, \;\; \s_i(y_j)=\l_{ij}'y_j, \;\; \tau_i(x_j)=\mu_{ij}x_j, \;\; \tau_i(y_j)=\mu_{ij}'y_j,$$
for some elements $\l_{ij}, \l_{ij}', \mu_{ij}$ and $\mu_{ij}'$ of the ring $D$. The elements $\L =(\l_{ij})$,  $\L' =(\l_{ij}')$,  $ M =(\mu_{ij})$ and  $M' =(\mu_{ij}')$ are called the {\em defining coefficients} of $A$. The $n$-tuples of endomorphisms $\s = (\s_1, \ldots , \s_n)$ and $\tau = (\tau_1, \ldots , \tau_n)$ are called the {\em defining endomorphisms} of $A$, and the $n$-tuple of  elements  $a=(a_1, \ldots , a_n)$ is called the {\em defining elements} of $A$. The GWA $A$ of rank $n$ is denoted by $A= D[x,y; \s , \tau, \L , \L', M, M']$ where $x=(x_1, \ldots , x_n)$ and $y=(y_1, \ldots , y_n)$.

We denote by $\s_i$ and $\tau_i$ the restrictions $\s_i|_D$ and $\tau_i|_D$, respectively.

An element $\L = (\l_{ij})$ (where $1\leq j <i\leq n$)  is called a {\em lower triangular half-matrix} with coefficients in $D$. The set of all such elements is denoted by $L_n(D)$. The next proposition describes GWAs of rank $n$ via generators and defining relations.

\begin{proposition}\label{28Jul}%\marginpar{28Jul}
Let $D$ be a ring, $\s = (\s_i)$ and $\tau = (\tau_i)$ be $n$-tuples of ring endomorphisms of $D$, $a=(a_i)\in D^n$, and $\L =(\l_{ij}), \L' =(\l_{ij}'),  M =(\mu_{ij}), M' =(\mu_{ij}')\in L_n(D)$ be such that the following conditions hold: For all $i=1, \ldots , n$ and $d\in D$,
%\marginpar{C1-3}
\begin{equation}\label{C1-3}
\tau_i\s_i(a_i) = a_i, \;\; a_id=\tau_i\s_i(d)a_i\;\; {\rm and}\;\; \s_i(a_i)d= \s_i\tau_i(d)\s_i(a_i);
\end{equation}
for  all $i>j$,
%\marginpar{C4-5}
\begin{equation}\label{C4-5}
a_i=\tau_i(\l_{ij})\mu_{ij}\s_j(a_i)=\tau_i(\l_{ij}') \mu_{ij}'\tau_j(a_i),
\end{equation}
%\marginpar{C6-7}
\begin{equation}\label{C6-7}
\s_i(a_i)=\s_i(\mu_{ij})\l_{ij}\s_j\s_i(a_i)=\s_i(\mu_{ij}')\l_{ij}'\tau_j\s_i(a_i);
\end{equation}
for all $i>j$ and $d\in D$,
%\marginpar{R1-2}
\begin{equation}\label{R1-2}
\l_{ij}\s_j\s_i(d)=\s_i\s_j(d)\l_{ij}\;\; {\rm and}\;\; \mu_{ij}\s_j\tau_i(d) = \tau_i\s_j(d) \mu_{ij},
\end{equation}
%\marginpar{R3-4}
\begin{equation}\label{R3-4}
\l_{ij}'\tau_j\s_i(d)=\s_i\tau_j(d)\l_{ij}'\;\; {\rm and}\;\; \mu_{ij}'\tau_j\tau_i(d) = \tau_i\tau_j(d) \mu_{ij}',
\end{equation}
%\marginpar{R5-6}
\begin{equation}\label{R5-6}
\s_i(a_j)=\l_{ij}'\tau_j(\l_{ij})a_j\;\; {\rm and}\;\; \tau_i(a_j)=\mu_{ij}'\tau_j(\mu_{ij})a_j,
\end{equation}
%\marginpar{R7-8}
\begin{equation}\label{R7-8}
\s_i\s_j(a_j)=\l_{ij}\s_j(\l_{ij}')\s_j(a_j)\;\; {\rm and}\;\; \tau_i\s_j(a_j)=\mu_{ij}\s_j(\mu_{ij}')\s_j(a_j).
\end{equation}
The GWA of rank $n$, $A= D[x,y; \s , \tau, a, \L , \L', M, M']$, is a ring generated by $D$, $x_1, \ldots , x_n$ and $y_1, \ldots , y_n$ subject to the defining relations: For all $i=1, \ldots , n$ and $d\in D$,
%\marginpar{DR1}
\begin{equation}\label{DR1}
x_id=\s_i(d)x_i, \;\; y_id= \tau_i(d) y_i, \;\; y_ix_i = a_i \;\; {\rm and}\;\; x_iy_i=\s_i(a_i);
\end{equation}
for all $i>j$,
%\marginpar{DR2}
\begin{equation}\label{DR2}
x_ix_j=\l_{ij}x_jx_i, \;\; x_iy_j=\l_{ij}'y_jx_i, \;\; y_ix_j=\mu_{ij}x_jy_i\;\; {\rm and}\;\; y_iy_j=\mu_{ij}'y_jy_i.
\end{equation}
\end{proposition}

{\it Proof}. The proof is routine and the result follows from the fact that a GWA of rank $n$ is a special type of the iterated GWA of rank $n$ and Theorem \ref{1Apr}.

Let $A= D[x,y; \s , \tau, a, \L , \L', M, M']$ be a GWA of rank $n$. Then the ring $A$ is generated by $D$, $x_1, \ldots , x_n$ and $y_1, \ldots , y_n$ subject to the defining relations (\ref{DR1}) and (\ref{DR2}) (by the very definition of $A$). The remaining equations (\ref{C1-3}), $\ldots$, (\ref{R7-8}) follow from  (\ref{abGWA}) and (\ref{GWADEF}) bearing in mind the iterated nature of the GWA $A$ and the definition of the endomorphisms $\s_1, \ldots ,\tau_n$. In more detail, the ring $A$ contains the obvious chain of iterated GWAs of rank $1,2,\ldots , n$: $A_1\subset A_2\subset\cdots \subset A_i\subset\cdots \subset A_n=A$ where $A_1=D[x_1, y_1; \ldots ], \ldots , A_i=A_{i-1}[x_i, y_i; \ldots ], \ldots , A=A_n=A_{n-1}[x_n, y_n; \ldots ]$. The relations (\ref{C1-3}) follow from (\ref{abGWA}). The relations (\ref{C4-5}) follow from the relations $a_id_{i-1}=\tau_i\s_i(d_{i-1})a_i$ in $A_i$ where $d_{i-1}=x_j$ and $ d_{i-1}=y_j$ $(i>j$), respectively (and using Theorem \ref{1Apr} for the iterated GWA $A_i=A_{i-1}[x_i, y_i; \ldots ]$):
\begin{eqnarray*}
d_{i-1}=x_j: & & a_ix_j=\tau_i\s_i(x_j)a_i=\tau_i(\l_{ij})\mu_{ij}\s_j(a_i)x_j, \\
 d_{i-1}=y_j: & & a_iy_j=\tau_i\s_i(y_j)a_i=\tau_i(\l_{ij}')\mu_{ij}'\tau_j(a_i)y_j.
\end{eqnarray*}
By Theorem \ref{1Apr}, we obtain the relations (\ref{C4-5}). Similarly, the relations (\ref{C6-7}) follow from the relations $\s_i(a_i) d_{i-1}=\s_i\tau_i(d_{i-1})\s_i(a_i)$ in $A_i$ where $d_{i-1}=x_j$ and $d_{i-1}=y_j$ $(i>j)$, respectively:
\begin{eqnarray*}
d_{i-1}=x_j: & & \s_i(a_i)x_j=\s_i\tau_i(x_j)\s_i(a_i)=\s_i(\mu_{ij})\l_{ij}\s_j\s_i(a_i)x_j, \\
 d_{i-1}=y_j: & & \s_i(a)_iy_j=\s_i\tau_i(y_j)\s_i(a_i)=\s_i(\mu_{ij}')\l_{ij}'\tau_j\s_i(a_i)y_j.
\end{eqnarray*}
By Theorem \ref{1Apr}, we obtain the relations (\ref{C6-7}).

The relations (\ref{R1-2}) follow from the relations $x_jd=\s_j(d) x_j$ (for all $d\in D$) in $A_i$ $(i>j)$ by applying the endomorphisms $\s_i$ and $\tau_i$ of the ring $A_{i-1}$, respectively (and using Theorem \ref{1Apr} for the iterated GWA $A_j=A_{j-1}[x_j, y_j; \ldots ]$:
\begin{eqnarray*}
\s_i: & & \l_{ij}\s_j\s_i(d) x_j=\s_i(x_jd)= \s_i(\s_j(d)x_j) = \s_i\s_j(d) \l_{ij}x_j,  \\
 \tau_i: & & \mu_{ij}\s_j\tau_i(d)x_j=\tau_i(x_jd)=\tau_i(\s_j(d)x_j) = \tau_i\s_j(d)\mu_{ij}x_j.
\end{eqnarray*}
By Theorem \ref{1Apr}, we obtain the relations (\ref{R1-2}).

Similarly, the relations (\ref{R3-4}) follow from the relations $y_jd=\tau_j(d)y_j$ (for all $d\in D$) in $A_i$ $(i>j)$ by applying the endomorphisms $\s_i$ and $\tau_i$ of the ring $A_{i-1}$, respectively:
\begin{eqnarray*}
\s_i: & & \l_{ij}'\tau_j\s_i(d) y_j=\s_i(y_jd)= \s_i(\tau_j(d)y_j) = \s_i\tau_j(d) \l_{ij}'y_j,  \\
\tau_i: & & \mu_{ij}'\tau_j\tau_i(d)y_j=\tau_i(y_jd)=\tau_i(\tau_j(d)y_j) = \tau_i\tau_j(d)\mu_{ij}'y_j.
\end{eqnarray*}
By Theorem \ref{1Apr}, we obtain (\ref{R3-4}).

The relations (\ref{R5-6}) follow from the relations $a_j=y_jx_j$ in $A_i$ $(i>j)$ by applying the endomorphisms $\s_i$ and $\tau_i$ of the ring $A_{i-1}$, respectively:
\begin{eqnarray*}
\s_i: & & \s_i(a_j) =\s_i(y_j)\s_i(x_j)=\l_{ij}'y_j\l_{ij}x_j=\l_{ij}'\tau_j(\l_{ij})a_j,   \\
\tau_i: & & \tau_i(a_j)=\tau_i(y_j)\tau_i(x_j)=\mu_{ij}'y_j\mu_{ij}x_j=\mu_{ij}'\tau_j(\mu_{ij})a_j.
\end{eqnarray*}
Similarly, the relations (\ref{R7-8}) are obtained from the relations $\s_j(a_j)=x_jy_j$ in $A_i$ $(i>j)$ by applying the endomorphisms $\s_i$ and $\tau_i$, respectively:
\begin{eqnarray*}
\s_i: & & \s_i\s_j(a_j) =\s_i(x_j)\s_i(y_j)=\l_{ij}x_j\l_{ij}'y_j=\l_{ij}\s_j(\l_{ij}')\s_j(a_j),   \\
\tau_i: & & \tau_i\s_j(a_j)=\tau_i(x_j)\tau_i(y_j)=\mu_{ij}x_j\mu_{ij}'y_j=\mu_{ij}\s_j(\mu_{ij}')\s_j(a_j).
\end{eqnarray*}
Conversely, suppose that the conditions of the theorem hold. Then using induction on $n$ we see that the ring $A$,  which is generated by $D$, $x_1, \ldots , x_n, y_1, \ldots , y_n$ and subject to the defining relations (\ref{DR1}) and (\ref{DR2}),  admits the chain of GWAs of rank $1,2,\ldots , n$, respectively:
$A_1\subset A_2\subset  \cdots \subset A_n=A$ (the relations (\ref{R1-2}), $\ldots , $ (\ref{R7-8}) are equivalent to the fact that the endomorphisms $\s_2, \ldots , \s_n$ of the ring $D$ can be lifted to the endomorphisms of the GWAs $A_1, \ldots , A_{n-1}$, respectively). $\Box$

{\em Example.} If $\s= (\s_1, \ldots , \s_n)\in \Aut (D)^n$ is an $n$-tuple of {\em commuting} automorphisms of the ring $D$, $\tau:= \s^{-1}= (\s_1^{-1}, \ldots , \s_n^{-1})$, $a=(a_1, \ldots , a_n)\in Z(D)$ and $\s_i(a_j) = a_j$ for all $i\neq j$; and $\l_{ij}=\l_{ij}'=\mu_{ij}=\mu_{ij}'=1$ for all $i>j$, then the GWA $A$ of rank $n$ is a {\em classical } GWA of rank $n$, that is $A=D[x,y;\s , a]$.

Recall that each normal, regular element $\alpha$ of a ring $D$ determines the  automorphism $\o_\alpha$ of $D$ by the rule $\alpha d = \o_\alpha (d) \alpha$ for all $d\in D$. The next proposition gives plenty of examples of GWAs of rank $n$.

\begin{proposition}\label{31Jul}%\marginpar{31Jul}
Let $D$ be a ring, $\th_1, \ldots , \th_n$  commuting automorphisms of the ring $D$, $\alpha_1, \ldots , \alpha_n$, $
\beta_1, \ldots , \beta_n$  regular, normal elements of $D$. Then $A= D[x,y; \s , \tau, a, \L , \L', M, M']$ be a GWA of rank $n$ where
\begin{eqnarray*}
\s_i &=& \th_i\o_{\beta_i}, \;\; \tau_i=\o_{\alpha_i}\th_i^{-1}, \;\; a_i=\alpha_i\beta_i, \\
\l_{ij} &=& \th_i(\beta_i)\th_i\th_j(\beta_j\beta_i^{-1})\th_j(\beta_j^{-1}), \;\; \l_{ij}'= \th_i(\beta_i\alpha_j)\cdot \th_j^{-1}\th_i(\beta_i^{-1})\alpha_j^{-1}, \\
 \mu_{ij} &=& \alpha_i\th_i^{-1}\th_j(\beta_j)\th_j(\alpha_i^{-1}\beta_j^{-1}), \;\; \mu_{ij}'= \alpha_i\th_i^{-1}(\alpha_j)\th_j^{-1}(\alpha_i^{-1})\alpha_j^{-1},
\end{eqnarray*}
provided $\l_{ij}, \l_{ij}', \mu_{ij}, \mu_{ij}\in D$.
\end{proposition}

{\it Proof}. It is routine to verify that the equalities (\ref{C1-3}), $\ldots $, (\ref{R7-8}) hold. Let us show that the first two equalities in (\ref{C1-3}) hold.
$$\tau_i\s_i(a_i)= \o_{\alpha_i}\th_i^{-1}\th_i\o_{\beta_i}(a_i)= \o_{\alpha_i\beta_i}(a_i)=\o_{a_i}(a_i)=a_i\;\; {\rm and} \;\; a_id=\o_{a_i}(d)d=\tau_i\s_i(d)d.$$
The third equality in (\ref{C1-3}) follows from the second by applying the automorphism $\s_i$ of the ring $D$.

Let us show that the first equality in (\ref{C4-5}) holds:
$$\tau_i(\l_{ij})\mu_{ij}\s_j(a_i)=\alpha_i(\beta_i\th_j(\beta_j\beta_i^{-1})\th_i^{-1}\th_j(\beta_j^{-1}))
\alpha_i^{-1}\cdot\alpha_i\th_i^{-1}\th_j(\beta_j)\th_j(\alpha_i^{-1}\beta_j^{-1})\cdot \th_j(\beta_j\alpha_i\beta_i\beta_j^{-1})=\alpha_i\beta_i=a_i.$$
In a similar fashion, the rest of the equalities are verified.  $\Box $

{\bf A $\Z^n$-grading of a GWA of rank $n$.} By Theorem \ref{1Apr}, every GWA of rank $n$, $A= D[x,y; \s , \tau, a, \L , \L', M, M']$, is a $\Z^n$-graded algebra $A=\oplus_{\alpha \in \Z^n} Dv_\alpha$ $(Dv_\alpha Dv_\beta\subseteq Dv_{\alpha +\beta}$ for all elements $\alpha , \beta \in \Z^n$) where for $\alpha = (\alpha_1, \ldots , \alpha_n)\in \Z^n$, $v_\alpha = v_{\alpha_1}(1)v_{\alpha_2}(2)\cdots v_{\alpha_n}(n)$ and
$v_{\alpha_i}(i):= \begin{cases}
x_i^{\alpha_i}& \text{if }\alpha_i\geq 0,\\
y_i^{-\alpha_i}& \text{if }\alpha_i<0.\\
\end{cases}$
Notice that the order in the product for $v_\alpha$ is important and, in general, cannot be changed. Moreover, the left $D$-module $Dv_\alpha$ is free of rank 1. For all elements $\alpha , \beta \in \Z^n$,
$$v_\alpha v_\beta = (\alpha , \beta ) v_{\alpha +\beta}$$
for some (explicit) elements $(\alpha, \beta ) \in D$. For all elements $\alpha \in \Z^n$ and $d\in D$, $v_\alpha d = \s^\alpha (d) v_\alpha$ where $\s^\alpha := \s(1, \alpha_1)\cdots \s (n, \alpha_n)$ and $\s (i, \alpha_i) := \begin{cases}
\s_i^{\alpha_i}& \text{if }\alpha_i\geq 0,\\
\tau_i^{-\alpha_i}& \text{if }\alpha_i<0.\\
\end{cases}$

%%%%%%%%%%%%%%%%%% SECTION 5  %%%%%%%%%%%%%%%%%%%%%%%%

\section{Diskew polynomial rings}\label{DPR}%\marginpar{DPR}

The aim of this section is to  show that  diskew polynomial rings are generalized Weyl algebras under a mild restriction (Theorem \ref{2Apr}), and to give proofs of simplicity criteria for them (Theorem \ref{Di31Oct}, Theorem \ref{QDi31Oct} and  Theorem \ref{PDi31Oct}).

{\bf Diskew polynomial rings}. Let $E= D[x,y; \s , \tau , b,\rho ]$ be a diskew polynomial ring.

{\it Remarks}.  1. If $\rho$ is a left normal element the two equalities in (\ref{wbd2}) can be written respectively as follows:  for all $d\in D$,
%\marginpar{wbd}
\begin{equation}\label{wbd}
(\s\tau (d) - \o_\rho \tau \s (d))\rho =0\;\; {\rm and}\;\; (\s \tau (d) - \o_b (d))b=0,
\end{equation}
where $\o_\rho$ and $\o_b$ are defined in (\ref{waLa}). So, in general, the elements $\rho$ and $b$ are not left regular in $D$.

2.  If, in addition, the elements $\rho$ and $b$ are left regular in $D$ then the conditions (\ref{wbd}) can be written as
%\marginpar{wbd1}
\begin{equation}\label{wbd1}
\s\tau = \o_\rho \tau \s = \o_b
\end{equation}
where $\o_\rho$ and $\o_b$ are ring endomorphisms of $D$ such that $\rho d = \o_\rho (d) \rho$  and $ bd = \o_b(d) b$ for all $d\in D$.

3. A particular case of diskew polynomial rings are the ambiskew polynomial rings: Let $\mF$ be a field, $A$ be an $\mF$-algebra, $\rho\in \mF\backslash \{ 0\}$ and let $\nu$ be a $\g$-normal element of $A$ for some $\mF$-automorphism $\g$ of $A$. Let $\alpha\in \Aut_\mF (A)$ be such that $\alpha \g = \g \alpha $  and let $\beta := \alpha^{-1}\g = \g \alpha^{-1}$, so that $\alpha \beta = \g = \beta \alpha$. Extend $\beta$ to an $\mF$-automorphism of $A[ y ; \alpha ] $ by setting $\beta (y) = \rho y$. By  {\cite[Exercise 2ZC]{GW-book}}, there is a $\beta$-derivation $\d$ of  $A[ y ; \alpha ] $ such that $\d (A)=0$ and $\d (y)=\nu$. The {\em ambiskew polynomial ring} $R(A, \alpha ,\nu , \rho )$ is the iterated polynomial ring  $A[ y ; \alpha ][x; \beta , \d ]$, see the paper of Jordan and Wells \cite{Jordan-Wells-2013} for details.

{\em Example.} Let $D=C[t_1, \ldots , t_n; \nu_1, \ldots , \nu_n]$  be a skew polynomial ring over a ring $C$, $\nu_1, \ldots , \nu_n$ commuting monomorphisms of the ring $C$. For all $\alpha\in \N^n$, $\nu^\alpha = \nu_1^{\alpha_1}\cdots \nu_n^{\alpha_n}$ is a monomorphism of $C$ and the element $t^\alpha= t_1^{\alpha_1}\cdots t_n^{\alpha_n}$ is a left normal, left regular element of the ring $D$ with $\o_\alpha = \nu^\alpha$. Each monomorphism $\nu_i$ of the ring $C$ can be extended to a monomorphism of the ring $D$ by the rule $\nu_i(t_j)=t_j$ for all $j$. Let $\rho$ be a central, $\nu$-invariant unit of $D$  (`$\nu$-invariant'  means $\nu_1(\rho ) = \cdots = \nu_n(\rho ) = \rho$). Let $\eta $ be a central, $\nu$-invariant, regular  element of $D$ and $b=\eta t^{\alpha + \beta}$. Let $\s = \nu^\alpha$ and $\tau = \nu^\beta$ where $\alpha , \beta \in \N^n$. Then the conditions (\ref{wbd1}) hold since $\o_\rho = \id_D$ and $\nu^{\alpha +\beta}= \s \tau = \tau \s = \o_b$. So, $E=D[x,y; \nu^\alpha , \nu^\beta, \eta t^{\alpha + \beta}, \rho ]$ is a diskew polynomial ring.

{\em Example.} The quantum plane $\L = K\langle p,q\, | \, pq= \l qp\rangle$ (over a field $K$ where $\l \in K^*$) is a skew polynomial ring $\L = K[q][p;\nu]$ where $\nu (q)= \l q$. Then $E=\L [x,y; \nu^\alpha , \nu^\beta, \eta t^{\alpha + \beta}, \rho ]$ is a diskew polynomial ring where $\eta , \rho \in K^*$ and $\alpha , \beta \in \N$ (see the previous example).

\begin{theorem}\label{5Apr}%\marginpar{5Apr}
The diskew polynomial ring  $E= D[x,y; \s ,\tau , b, \rho]$ is an iterated skew polynomial ring $E=D[y; \tau][x; \s ,\der ]$ where $\s (y) = \rho y$, $\der (D)=0$ and $\der (y) =b$. It is a free left $D$-module $E= \oplus_{i, j\in \N} Dy^ix^j$ and the element $x$ is a left regular element.
\end{theorem}

{\it Proof}. Using (\ref{wbd2}), we have to show   that $\s$ and $\der$ respect the defining relation $yd-\tau (d) y = 0$ (for $d\in D$) of the  skew polynomial ring $D[y; \tau]$:
\begin{eqnarray*}
\s (yd-\tau (d) y) &=& \rho y\s (d) -\s \tau (d)\rho y=(\rho\tau\s (d) - \s\tau (d) \rho)y=0,\;\;  {\rm by}\;\; (\ref{wbd2}), \\
\der (yd-\tau (d) y) &=&  bd-\s\tau (d)b=0, {\rm by}\;\; (\ref{wbd2}). \;\; \Box
 \end{eqnarray*}

{\bf The $(x,y)$-symmetry of diskew polynomial rings.}  If $\rho\in D$ is a {\em unit} then the equality $xy-\rho yx=b$ is equivalent to the equality $yx-\rho^{-1}xy=-\rho^{-1}b$. Therefore,
%\marginpar{Exysym}
\begin{equation}\label{Exysym}
E=D[x,y;\s , \tau , b,\rho]=D[y, x; \tau , \s , -\rho^{-1}b, \rho^{-1}].
\end{equation}
When we say the $(x,y)$-symmetry of diskew polynomial rings we mean the second equality in (\ref{Exysym}). So, when $\rho$ is a unit the properties of diskew polynomial rings are more symmetrical.

{\bf  Proof of Theorem \ref{2Apr}.}  By the assumption, $\rho$ is a unit. Then  the elements $x$ and $y$  are  left regular element of $E$, by Theorem \ref{5Apr} and (\ref{Exysym}). Therefore, the subring of $E$ generated by $D$ and the left regular element $h=yx$ is the skew polynomial ring $\CD := D[h;  \tau \s]$.  Since  $\rho$ is a unit,  the first condition in (\ref{wbd}) can be written as $\s \tau = \o_\rho \tau \s$. This equality is used in the proof.

(i) $\s, \tau \in \End (\CD )$: We have to show that $\s$ and $\tau$ respect the defining relations, $hd = \tau \s (d) h$ ($d\in D$), of the ring $\CD$:
$$\s (hd-\tau \s (d)h) =  (\rho h +b) \s (d) - \s \tau \s (d) (\rho h +b)=(\o_\rho \tau \s \s (d) - \s \tau \s (d)) \rho h + (\o_b \s (d) - \s \tau \s (d))b\stackrel{(\ref{wbd})}{=}0.$$
\begin{eqnarray*}
\tau (hd-\tau \s (d)h) &= & \tau (\rho^{-1}) (h-\tau (b)) \tau (d) - \tau^2\s (d) \tau (\rho^{-1}) (h-\tau (b)) = \tau (\rho^{-1} \s \tau (d) - \tau \s (d) \rho^{-1} ) h\\
 &&  - \tau (\rho^{-1} bd-\tau \s (d) \rho^{-1} b) =\tau (\rho^{-1} (\s \tau (d) -\rho \tau \s (d) \rho^{-1} ))h-\tau (\rho^{-1} (\o_b(d) - \rho \tau \s (d) \rho^{-1})b)\\
 &= & \tau (\rho^{-1} (\s \tau (d) - \o_\rho \tau \s (d)))h-\tau (\rho^{-1} (\o_b(d)-\o_\rho \tau \s (d) )b)\\
 &=& -\tau (\rho^{-1}(\o_b(d) -\s\tau (d))b)=0,\;\;\; ({\rm by}\;\; (\ref{wbd})).
\end{eqnarray*}

(ii) {\em The conditions in (\ref{abGWA}) hold}: Recall that $a=h$. The first condition in (\ref{abGWA}), $\tau \s (a) = a$,  holds:
$$ \tau \s (h) = \tau ( \rho h+b) = \tau (\rho ) \tau (\rho^{-1}) (h-\tau (b)) +\tau (b) = h.$$
The second condition in (\ref{abGWA}), $ ad= \tau \s (d) a$,  holds since it means that $hd = \tau \s (d) h$. The third condition in (\ref{abGWA}), $\s (a) d = \s \tau (d) \s (a)$,  holds:
$$ \s (h) d = (\rho h+b) d = \o_\rho \tau \s (d) \rho \cdot h +\o_b(d) b  \stackrel{(\ref{wbd})}{=} \s \tau (d) \rho \cdot h +\s \tau (d)b= \s\tau (d) (\rho h+b) = \s\tau (d) \s (h).$$
(iii) {\em The defining relations (\ref{GWADEF}) of the GWA  hold}: The first two equalities  in (\ref{GWADEF}) are given. Now, $yx=h$ and $xy=\rho yx +b=\s (h)$.

(iv) $E= \CD [ x,y; \s , \tau , a=h]$ {\em is a  GWA}: By (i)-(iii), the ring $E$ is an epimorphic image of the GWA $A:= \CD [ x,y; \s , \tau , a=h]$. Since $E=\bigoplus_{i\geq 1}\CD y^i\oplus\bigoplus_{i\geq 1} \CD x^i$ (as $\rho$ is a unit) and $A= \oplus_{i\geq 1} \CD y^i\oplus\bigoplus_{i\geq 1} \CD x^i$, the epimorphic image is, in fact, an isomorphic image, i.e., $E\simeq A$.

(v) Since $\rho\in D$ is a unit, by (\ref{wbd}), $\s \tau = \o_\rho \tau \s$, and so $\s\tau (h) = \o_\rho (h) = \rho h \rho^{-1} = \rho \tau \s (\rho^{-1}) h$.

(vi) Since $\rho\in D$ is a unit, $\s\tau = \o_\rho \tau \s$, by (\ref{wbd}). This is an equality of endomorphisms of the ring $D$. This is also an equality of endomorphisms of the ring $\CD$: $\s\tau (h) = \o_\rho (h)$ and $\tau \s (h) = h$ imply $\s \tau (h) = \o_\rho \tau \s (h)$.  $\Box $

\begin{corollary}\label{a2Apr}%\marginpar{a2Apr}
Let $E=D[x,y;\s , \tau , b,\rho]$ be a diskew polynomial ring. Suppose that $\rho$ is a unit in $D$. Then $E= \CD [ y,x; \tau , \s , h':= \s (h)=\rho h+b]$ is a GWA with base ring $\CD := D[h; \tau \s ] = D[ h', \s\tau ]$ which is a skew polynomial ring, $\s$ and $\tau$ are ring endomorphism of $\CD$ that are extensions of the ring endomorphisms $\s$ and $\tau$ of $D$, respectively,  defined in Theorem \ref{2Apr}; $\tau (h') = \rho^{-1} (h'-b)$ and $\s (h')= \s (\rho ) h'+\s (b)$. In particular, $\s\tau (h') = h'$ and $\tau \s (h') = \o_{\rho^{-1}}(h') = \rho^{-1} \s\tau (\rho ) h'$.
\end{corollary}

{\it Proof}. By Theorem \ref{2Apr}, $E=\CD [x,y; \s , \tau , h]$ where $\CD = D[h;\tau \s ]$, $\s (h) = \rho h+b$ and $ \tau (h)=\tau (\rho^{-1})(h-\tau (b))$. By the $(x,y)$-symmetry for GWAs, $\CD [ x,y; \s , \tau , h]= \CD [ y,x; \tau , \s , \s (h) = \rho h+b]$. By the $(x,y)$-symmetry for DPRs, $E= D[y,x; \tau , \s , -\rho^{-1} b, \rho^{-1}]$. Then, by Theorem \ref{a2Apr}, $E= D[y,x; \tau , \s , -\rho^{-1} b, \rho^{-1}]= \CD'[y,x; \tau , \s , \s (h)]$ where $\CD' = D[\s (h) , \s\tau ]=D[h, \tau \s ]=\CD$. Now,
 $\tau (h') = \tau \s (h)\stackrel{{\rm Th}\, \ref{2Apr}}{=} h= \rho^{-1} (\rho h+b-b)=\rho^{-1}(\s (h)-b) = \rho^{-1} (h'-b)$ and $\s (h')= \s (\rho h+b) = \s (\rho ) \s (h)+\s (b) = \s (\rho ) h'+\s (b)$. Furthermore,  $\s \tau (h') = \s \tau \s (h) = \s (h) = h'$ (since $ \tau \s (h) = h$, by Theorem \ref{2Apr}) and $\tau \s (h') = \o_{\rho^{-1}}\s\tau (h') = \o_{\rho^{-1}}\s\tau \s (h) = \o_{\rho^{-1}}\s (h) = \o_{\rho^{-1}}(h')= \rho^{-1} h' \rho = \rho^{-1} \s\tau (\rho ) h'$.
 $\Box $

 Under mild conditions, the next corollary produces a series of diskew polynomial rings from a given one. The construction is based on Corollary \ref{a2Apr}.

\begin{corollary}\label{a7AApr}%\marginpar{a7AApr}
Let $E=D[x,y;\s , \tau , b,\rho]$ be a diskew polynomial ring such that $\rho$ is a unit in $D$.
\begin{enumerate}
\item Let $\CD_1= D[h_1; \nu := \tau \s ]$ be a skew polynomial ring where $h_1$ is a variable. The endomorphisms  $\s$ and $\tau$ of the ring $D$ can be extended to  endomorphisms of the ring  $\CD$ by the rule  $\s (h_1)= \rho  h_1+b$ and  $\tau (h_1) =\tau ( \rho^{-1}) (h_1-\tau (b))$. In particular, $\tau \s (h_1) = h_1$,  $\s\tau  (h_1) = \o_{\rho}(h_1) = \rho \tau \s (\rho^{-1} ) h_1$ and $\s\tau=\o_\rho \tau \s$ in $\CD$. The endomorphisms $\s$ and $\tau$ of $\CD$ satisfy (\ref{wbd2}) iff the element $\rho^{-1}b$ is $\tau\s$-invariant.
\item Let $\CD_n= D[h_1, \ldots , h_n; \tau\s, \ldots , \tau \s ]$ be a skew polynomial ring such that the element $\rho^{-1}b$ is $\tau\s$-invariant. Then $E_n= \CD_n[x,y;\s , \tau , b,\rho ]$ is  a diskew polynomial ring where $\s (h_i) = \rho h_i+b$ and $\tau (h_i)= \tau(\rho^{-1})(h_i-\tau (b))$ for $i=1, \ldots , n$, and, by Theorem \ref{2Apr}, the ring $E_n=\CD_{n+1} [x,y; \s ,\tau, h_{n+1}]$ is a GWA.
\end{enumerate}
\end{corollary}

{\it Proof}. 1. Repeating word for word the proof of Theorem \ref{2Apr}, we obtain all the statements in statement 1 but the last sentence. The endomorphisms $\s$ and $\tau$ satisfy (\ref{wbd2}) iff $\s\tau (h_1) b=bh_1$. Since the element $h_1$ is left regular, this equality holds iff the element $\rho^{-1}b$ is $\tau \s$-invariant since $ \s\tau (h_1) b= \rho \tau \s (\rho^{-1}) h_1b= \rho \tau \s (\rho^{-1}) \tau \s (b) h_1= \rho \tau \s (\rho^{-1} b) h_1$.

2. Using repeatedly  statement 1, we obtain the iterated skew polynomial ring $\CD_n= D [h_1; \tau\s ]$ $\cdots [h_n; \tau \s ]=D[h_1, \ldots , h_n; \tau\s, \ldots , \tau \s ]$ (since $ h_{i+1}h_i= \tau\s (h_i) h_{i+1}=h_ih_{i+1}$) such that the extension of the endomorphisms $\s$ and $\tau$ from $D$ to $\CD_n$ ( as in statement 2) satisfy (\ref{wbd2}). Now,  statement 2 is obvious. $\Box$

By Theorem \ref{2Apr}, if $\rho\in D$ is a unit then the ring endomorphisms $\s$ and $\tau$ of $D$ can be extended to the ring $\CD = D[h;  \tau\s ]$ by the rule $\s (h) = \rho h +b$ and $\tau (h) = \tau (\rho^{-1}) (h-\tau (b))$. By induction on $i\geq 1$, we have the equalities where
%\marginpar{sih1}
\begin{equation}\label{sih1}
\s^i(h) = a_ih+b_i, \;\; a_i = \s^{i-1} (\rho ) \cdots \s (\rho ) \rho \;\; {\rm and}\;\; b_i=\sum_{j=1}^{i-1} \s^j(a_{i-j})\s^{j-1} (b) +\s^{i-1}(b),
\end{equation}
%\marginpar{sih2}
\begin{equation}\label{sih2}
\tau^i(h) = a_i'h+b_i', \;\; a_i' = \tau^{i-1} (\rho^{-1} ) \cdots \tau^2 (\rho^{-1} ) \tau (\rho^{-1} ) \;\; {\rm and}\;\; b_i'=-\sum_{j=1}^{i-1} \tau^j(a_{i-j}')\tau^j (\rho^{-1}b) -\tau^i(\rho^{-1}b).
\end{equation}
In particular, for all $i\geq 1$,
%\marginpar{sih3}
\begin{equation}\label{sih3}
a_{i+1}= \s (a_i) \rho \;\; {\rm and}\;\;b_{i+1}= \s (a_i) b+\s (b_i),
\end{equation}
%\marginpar{sih4}
\begin{equation}\label{sih4}
a_{i+1}'= \tau (a_i') \tau (\rho^{-1}) \;\; {\rm and}\;\; b_{i+1}'= -\tau (a_i')\tau (\rho^{-1} b)+\tau (b_i').
\end{equation}
For example, (\ref{sih3}) follows from $a_{i+1}h+b_{i+1}= \s ( \s^i(h))= \s ( a_ih+b_i) = \s (a_i) (\rho h +b) + \s (b_i) = \s (a_i) \rho h + \s (a_i) b + \s (b_i)$. For two elements $s$ and $t$ of a ring, $[s,t]:= st-ts$ is its {\em commutator}.

Suppose that $\rho$ is a {\em unit}. Then, by (\ref{wbd}), $\s\tau = \o_\rho\nu$ where $\nu = \tau \s$, or, equivalently, $\o_{\rho^{-1}}\s\tau = \nu$. Let $\beta :=\rho^{-1}b$. It follows   that for all $d\in D$,
%\marginpar{bdnb}
\begin{equation}\label{bdnb}
\beta d= \nu (d) \beta\;\; {\rm and}\;\; (h+\beta ) d= \nu (d)(h+\beta )
\end{equation}
$(\beta d = \rho^{-1} bd\stackrel{(\ref{wbd2})}{=}\rho^{-1} \s\tau (d) b = \o_{\rho^{-1}}\s\tau (d) \beta = \nu (d)\beta )$.
 If, in addition, we assume that the element $b$ is a {\em left regular} element $D$.  Then, by (\ref{wbd}), $\s\tau = \o_\rho \nu = \o_b$ and the element $\beta\in D$ is also  left regular in $D$. By (\ref{bdnb}), $\beta \beta = \nu (\beta ) \beta$ and $ (h+\beta ) \beta = \nu (\beta) (h+\beta )$. Hence,
 %\marginpar{bdnb1}
\begin{equation}\label{bdnb1}
\nu (\beta )= \beta , \;\; h\beta = \beta h ,   \;\;  \s (h^i) = \rho_i^\nu\sum_{j=0}^i\beta_{ij} h^j \;\; (i\geq 1) \; {\rm where}\;  \rho_i^\nu :=\rho \nu(\rho) \cdots \nu^{i-1}(\rho ), \; \beta_{ij}={i\choose j}\beta^{i-j}.
\end{equation}
In more detail, $\s (h^i) = (\rho (h+\beta ))^i= \rho_i^\nu (h+\beta )^i =\rho_i^\nu\sum_{j=0}^i\beta_{ij} h^j$. For all natural numbers $i$ and $j$ such that $1\leq j\leq i$, $(\rho^\nu_j)^{-1} \rho_i^\nu = \nu^j (\rho^\nu_{i-j})$.

{\bf Simplicity criterion for diskew polynomial rings where $\rho$ is a unit.} Theorem \ref{Di31Oct} is a simplicity criterion for a diskew polynomial ring $E=D[x,y;\s , \tau , b, \rho ]$ where $\rho$ is a unit and $\tau \s$ is an epimorphism. By Theorem \ref{2Apr}, the ring $E$ is a GWA that satisfies the assumptions of Theorem \ref{B5Apr},  a simplicity criterion for GWAs, and Theorem \ref{Di31Oct} is rather a straightforward corollary of Theorem \ref{B5Apr}.

{\bf Proof of Theorem \ref{Di31Oct}}. $(1\Leftrightarrow 2)$ Since $\rho$ is a unit, $E=\CD [x,y;\s , \tau , a=h ]$ is a GWA (by Theorem \ref{2Apr}) where $\CD = D[h; \tau \s ]$. Since $\tau \s$ is an epimorphism, the element $a=h$ is a {\em normal} element of the ring $\CD$. By Corollary \ref{a2Apr}, $E= \CD [ y,x;\tau , \s , \s (a)]$ and $\CD = D[\s (h) ; \s\tau ]$. By Theorem \ref{2Apr}, $\s \tau =\o_\rho \tau \s$ in $\CD$. Hence, $\s \tau$ is an epimorphism of the ring $\CD$, and so the element $\s (a) = \s (h)$ is a normal element of $\CD$.  This means that the GWA $E= \CD [ x,y;\s , \tau , a=h]$ satisfies the assumptions of Theorem \ref{B5Apr}. In particular, the ring $E$ is a simple ring iff the conditions (a)-(d) of Theorem \ref{B5Apr} hold for the GWA $E$. We aim to show that the conditions (a)-(d) of statement 2 of the theorem  are equivalent to the conditions (a)-(d) of Theorem \ref{B5Apr}.

$(a) \Leftrightarrow (a):$ Since $\tau\s$ and $\s\tau$ are epimorphisms of the ring $\CD = D[h; \tau \s]= D[\s (h) ; \s\tau ]$, the conditions that the elements $a$ and $\s (a)$ are regular in $\CD$ are equivalent to the conditions that $\tau \s$ and $\s \tau$ are automorphisms of $\CD$ or, equivalently, $\s$ and $\tau$ are automorphisms of $D$ (since $\s (h) = \rho h+b$ and $\tau (h) = \tau (\rho^{-1})(h-\tau (b))$, by Theorem \ref{2Apr}).

$(a,d) \Leftrightarrow (a,d):$  Since $\CD = D[ h; \tau \s ]$, the condition $\CD h +\CD \s^i (h)=\CD$ (where $i\geq 1$) holds iff $\CD b_i= \CD$ (by (\ref{sih1})) iff $Db_i=D$ iff $b_i^*b_i=1$ for some element $b_i^*\in D$.

Since $\tau\s$ is an automorphism of the ring $D$ (as $\s$ and $ \tau$  are automorphisms), the element $h$ of $\CD$ is a normal, regular element. Hence, so is the element $\s^i(h)$ of $\CD$ (since $\s$  is an automorphism of the ring $\CD$). Then the condition $\CD h +\CD \s^i (h) = \CD$ can be rewritten as the condition $h\CD +\s^i (h) \CD = \CD$ which is equivalent to the equality $b_iD= D$, i.e., $b_ib_i^o=1$ for some element $b_i^0\in D$. Therefore, the condition $\CD h +\CD \s^i (h) = \CD$ is equivalent to the conditions that the element $b_i$ of $D$ is a unit (since $(a) \Leftrightarrow (a)$).

$(a,c, d) \Leftrightarrow (a,d):$  In view of the equivalence `$(a, d) \Leftrightarrow (a,d)$', it suffices to show that the implication $(c) \Leftarrow (a,d)$ holds, i.e., that the automorphism $\s^n$ of $\CD$ is not inner for all $n\geq 1$. Suppose that this is not the case for some $n\geq 1$, that is $ud= \s^n(d) u$ for all $d\in \CD$ where $u=\sum_{i=0}^mu_ih^i$ is a unit of $\CD$ and $u_i\in D$. Then necessarily the element  $u_0$ is a unit of $D$ (since $\CD / (h) \simeq D$). Taking the equality $uh= \s^n (h) u= (a_nh+b_n) u$ modulo the ideal $(h)$ of the ring $\CD$ we have that $b_nu_0=0$ in $D$ (since $\CD / (h) \simeq D$), hence $u_0=0$ (since $b_n$ is a unit of $D$), a contradiction.

$(a,b,c,d) \Leftrightarrow (a,b,c,d):$  We know already that $(a,c,d) \Leftrightarrow (a,d)$.  Suppose that the ring  $\CD$ is a $\s$-simple ring. Then necessarily {\em the ring $D$ is a $\s$-simple ring}:  If $\gp$ is a nonzero, $\s$-invariant ideal of the ring $D$ then $\CD \gp \CD$ is a nonzero, $\s$-invariant  ideal of the ring $\CD$. Therefore, $\CD \gp \CD = \CD$. By taking this equality modulo the ideal $(h)$ of the ring $\CD$, we have the equality $\gp = D$, as required.

Suppose that $D$ is a $\s$-simple ring.
Let $I$ be a nonzero, $\s$-invariant ideal of the ring $\CD$. Let $p= \sum_{i=0}^n d_ih^i$ (where $d_i\in D$) be a nonzero element of $I$ of {\em least possible degree} $n$ with respect to $h$. In particular, $d_n\neq 0$ and $n\geq 1$ (since $D$ is a $\s$-simple ring). Then,
\begin{eqnarray*}
I\supseteq  \sum_{i\geq 0} D\s^i (p) D &=& \sum_{i\geq 0} D\s^i(d_n) a_ih^nD+\cdots = \bigg( \sum_{i\geq 0} D\s^i (d_n)a_i\nu^n (D)\bigg) h^n+\cdots\\
  &=& \bigg( \sum_{i\geq 0} D\s^i (d_n)D\bigg) h^n+\cdots = Dh^n+\cdots
\end{eqnarray*}
since the elements $a_i$ are units of $D$ (see (\ref{sih1})), $\nu$ is an automorphism of $D$ and $\sum_{i\geq 0} D\s^i (d_n)D$ is a nonzero, $\s$-invariant ideal of the ring $D$ which is equal to $D$ (by the $\s$-simplicity of $D$). Therefore, without loss of generality we may assume that $b_n=1$, i.e., $p=h^n+\sum_{i=0}^{n-1}d_ih^i$.

(i) {\em For all elements $d\in D$, $pd= \nu^n (d) p$, i.e., $d_id= \nu^{n-i}(d) d_i$ for} $i=0,1,\ldots , n-1$: The element of the ideal $I$,
$$ \nu^n (d) p-pd= \sum_{i=0}^{n-1}(\nu^n(d) d_i-d_i\nu^i(d))h^i, $$
has degree $<n$, hence it is equal to zero, by the minimality of $n$, and the statement (i) follows (since $\nu$ is an automorphism of $D$, the element $h$ is regular in $\CD$).

(ii) $\s (p) = \rho_n^\nu  p\, $: The element of the ideal $I$,
$$ \s (p) - \rho_n^\nu p = (\rho_n^\nu h^n+\cdots ) - (\rho_n^\nu h^n+\cdots )$$ has degree $<n$, hence  it is equal to zero, by the minimality of $n$, and the statement (ii) follows.

(iii) $[h,p]=0$: For $i=0, 1, \ldots , n-1$, $d_i\beta = \nu^{n-i}(\beta ) d_i = \beta d_i$ since $\nu (\beta ) = \beta$, by (\ref{bdnb1}). By (\ref{bdnb1}), $\beta h = h\beta$. Hence,  $p\beta = \beta p$. The element $\beta = \rho^{-1} b$ is a unit since $\rho$ and $b$ are so (see the statement (d)). Then, by (\ref{bdnb}), $d_i\beta = \beta d_i = \nu (d_i) \beta$, and so $d_i = \nu (d_i)$. In summary,
%\marginpar{pbbp}
\begin{equation}\label{pbbp}
p\beta = \beta p, \;\; \nu (d_i) = d_i\;\; {\rm and}\;\; d_i\beta = \beta d_i\;\; {\rm for}\;\; i=0,1, \ldots , n-1.
\end{equation}
Now, $hp=(h^n+\sum_{i=0}^{n-1}\nu (d_i) h^i) h= (h^n+\sum_{i=0}^{n-1}d_i h^i) h=ph$.

The statements (i)  and (ii) mean that the element $p$ is a regular, normal element of the ring $\CD$. Hence, $\CD p \CD = \CD p = p\CD$. The condition (ii) implies that the ideal $\CD p$ of $\CD$ is a proper, $\s$-invariant ideal of $\CD$. Now, it is obvious that the ring $\CD$ is a $\s$-simple ring iff the conditions (b) and (c) hold.

$(1\Leftrightarrow 3)$ In view of Corollary \ref{a2Apr},  this equivalence follows from the equivalence  $(1\Leftrightarrow 2)$  by the $(x,y)$-symmetry, the fact that the conditions that $\rho$ is a unit and $\tau \s$ is an epimorphism are equivalent to the conditions that $\rho^{-1}$ is a unit and $\s\tau$ is an epimorphism of the ring $D$ since $\s\tau = \o_\rho \tau \s$, by (\ref{wbd}).  $\Box $

We keep the notation and assumptions of Theorem \ref{Di31Oct}. Let  $p=h^n+\sum_{i=0}^{n-1}d_ih^i\in \CD$. Since $\rho$ is a unit,   then, by (\ref{bdnb1}), {\em the condition that $\s (p)=\rho_n^\nu p$ is equivalent to the equalities} (since $\s (p)=\rho_n^\nu h^n+\cdots $)
%\marginpar{hpdn2}
\begin{equation}\label{hpdn2}
\rho^\nu_n\beta_{nj}+\sum_{j\leq i\leq n-1}\s(d_i) \rho_i^\nu \beta_{ij}= \rho_n^\nu d_j \;  \; {\rm for}\; \; j=0,1, \ldots , n-1.
\end{equation}
For each $j$, by multiplying the equality (\ref{hpdn2}) by the unit $(\rho^\nu_j)^{-1}$ on the left we have the equality
$$\rho^\nu_{n-j}d_j-\s (d_j) = \rho^\nu_{n-j}\beta_{nj}+\sum_{j< i\leq n-1}\s(d_i) \rho_{i-j}^\nu \beta_{ij} \;  \; {\rm for}\; \; j=0,1, \ldots , n-1,$$
since $\rho_n^\nu d_j(\rho_j^\nu)^{-1} = \rho_n^\nu \nu^{n-j}(\rho_j^\nu )^{-1} d_j= \rho^\nu_{n-j}d_j$ and $\rho^\nu_i\beta_{ij}(\rho^\nu_j)^{-1} = \rho_i^\nu (\o_{\rho^{-1}}\o_b)^{i-j} ((\rho_j^\nu)^{-1})\beta_{ij} = \rho_i^\nu \nu^{i-j}(\rho_j^\nu )^{-1}\beta_{ij}= \rho^\nu_{i-j}\beta_{ij}$. Since $b^i=\rho^\nu_i\beta^i$ for all $i\geq 1$, the equalities above can be written as follows

%\marginpar{hpdn7}
\begin{equation}\label{hpdn7}
\rho^\nu_{n-j}d_j-\s (d_j) = {n\choose j}b^{n-j}+\sum_{j<i\leq n-1}{i\choose j}\s(d_i) b^{i-j}\;  \; {\rm for}\; \; j=0,1, \ldots , n-1.
\end{equation}

In particular, for $j=n-1$ and $j=0$ we have, respectively, the equalities
%\marginpar{hpdn6}
\begin{equation}\label{hpdn6}
\rho d_{n-1} - \s (d_{n-1})= nb,
\end{equation}
%\marginpar{hpdn4}
\begin{equation}\label{hpdn4}
\rho_n^\nu d_0-\s (d_0)=  b^n+ \sum_{1\leq  i\leq n-1} \s (d_i) b^i.
\end{equation}
Let $p=h^n+\sum_{i=0}^{n-1}d_ih^i\in \CD$ be as in Theorem \ref{Di31Oct}. Notice that the element $\beta$ is a unit, $\nu = \o_{\rho^{-1}b}=\o_{\beta}$ and $p\beta = \beta p$ (see (\ref{pbbp}). Then, $p=\o_{\beta}(p) = \nu (p) = \tau \s (p) = \tau (\rho_n^\nu p)$, and so
%\marginpar{pbbp1}
\begin{equation}\label{pbbp1}
\tau (p) = \tau(\rho_n^\nu)^{-1}p.
\end{equation}

{\bf Rings with enough normal elements.} We say that a ring has {\em enough normal elements} if each nonzero ideal contains a normal element. All commutative rings have enough normal elements. In a similar way, a ring that has {\em enough left/right normal elements} is defined. The next corollary provides examples of DPRs/GWAs that have enough regular normal elements.

\begin{corollary}\label{a8Nov}%\marginpar{a8Nov}
Let $E=D[x,y;\s , \tau , b, \rho ]$ be a diskew polynomial ring such that $\rho$ is a unit in $D$ and $\nu = \tau \s$ is an epimorphism. Suppose that the ring $D$ is  $\s$-simple (resp.,  $\tau$-simple); $\s$ (resp.,  $\tau$)  is an automorphism of $D$  and the elements $b_i$ (resp.,  $b_i'$) are units in $D$ for all $i\geq 1$. Then
 \begin{enumerate}
\item  every proper ideal of $E$ contains an element $p$ (resp., $p'$) that satisfies the conditions 2(c) (resp.,  3(c))  of Theorem \ref{Di31Oct}. In particular, the element $p$ (resp., $p'$) is a regular, normal element. The element $p$  (resp., $p'$) is unique provided its $h$-degree (resp., $h'$-degree) is the least possible.
\item The ring $E$ has enough regular normal elements.
\item The multiplicative monoid $\CP$ (resp., $\CP'$) generated by all the elements $p$ (resp., $p'$) is a regular  Ore set in $E$ such that the ring $\CP^{-1}E$ (resp., $\CP'^{-1}E$) is a simple ring.
  %  \item Every nonzero prime ideal of the ring $E$ contains a prime ideal generated by a left normal element %$p$ (resp., $p'$).
\end{enumerate}
\end{corollary}

{\it Proof}. The assumptions of the corollary are precisely the conditions (a), (b), (d) in statements 2 and 3 of Theorem \ref{Di31Oct}.  Now, statement 1 of  the corollary follows, see the proof of $(a,b,c,d)\Leftrightarrow (a,b,c,d)$ of Theorem \ref{Di31Oct}. Statements 2 and 3 follow from directly from statement 1.  $\Box$

The next theorems shed  light on the elements $p$ and $p'$ in Theorem \ref{Di31Oct}. Theorem \ref{A6Apr} describes the element $p$ in Theorem \ref{Di31Oct} where $n=1$.

\begin{theorem}\label{A6Apr}%\marginpar{A6Apr}
Let $E= D[x,y; \s , \tau , b, \rho ]$ be a diskew polynomial ring such that $\rho$ is a unit and $\CD = D[h; \nu = \tau \s ]$ where $h=yx$.  The following statements are equivalent.
\begin{enumerate}
\item There exists an element $C= h+\alpha \in \CD $, where $\alpha \in D$, such that $Cd=\nu (d) C$ for all elements $d\in D$ and $\s (C) = \rho C$.
\item There is an element $\alpha \in D$ such that $\rho \alpha - \s (\alpha ) = b$  and $\alpha d= \nu (d) \alpha$ for all elements $d\in D$.

If one of the equivalent conditions holds then $[h,C]=(\nu (\alpha ) - \alpha ) C$ and
\begin{enumerate}
\item  $C=\rho^{-1} (xy+\s (\alpha ))$, $xC= \rho Cx$ and $yC= \tau (\rho^{-1}) (C+\nu (\alpha ) -\alpha ) y$.
\item $E\simeq D[C; \nu ] [x,y; \s , \tau , a:= C-\alpha ]$ is a GWA where $\s (C) = \rho  C$ and $\tau (C) = \tau (\rho^{-1}) (C+\nu (\alpha ) -\alpha ) $. Furthermore, $\tau \s (C) = C+\nu (\alpha ) - \alpha$ and $\s \tau (C)=\s \tau (\rho^{-1})(\rho C+\s (\nu (\alpha ) - \alpha ))$.
\end{enumerate}
\end{enumerate}
\end{theorem}

{\it Proof}. $(1\Leftrightarrow 2)$ The equality $Cd= \nu (d) C$ is equivalent to the equality $\alpha d= \nu (d) d$. The equality $\s (C)= \rho C$, i.e., $\rho h+b+\s (\alpha ) = \rho h +\rho \alpha$ is equivalent to the equality $\rho\alpha - \s (\alpha ) = b$. Now, $[h,C]= (h+\nu (\alpha ) - C) h= (\nu (\alpha ) - \alpha ) h= (\nu (\alpha ) - \alpha )C- (\nu (\alpha ) - \alpha )\alpha = (\nu (\alpha ) - \alpha )C$ since $\alpha \cdot \alpha = \nu (\alpha ) \cdot \alpha$.

(a) $C= \rho^{-1} (xy+\s (\alpha ))$: $C= yx+\alpha = \rho^{-1}(xy-b) +\alpha = \rho^{-1} (xy+\s (\alpha ))$ since $ \rho \alpha - \s (\alpha ) = b$.

 $xC = x(yx+\alpha ) = (xy+\s (\alpha )) x= \rho Cx$ and $yC= y\rho^{-1} (xy+\s (\alpha )) = \tau (\rho^{-1})(yx+\tau \s (\alpha )) y= \tau (\rho^{-1})(C+\nu (\alpha ) - \alpha )y$.

(b) By Theorem \ref{2Apr}, $E= \CD [x,y; \s , \tau , h]$ is a GWA. Now, the statement (b) follows from the statement (a). In particular, $\tau \s (C) = \tau (\rho ) \tau (C) = \tau (\rho ) \tau (\rho^{-1}) (C+\nu (\alpha ) - \alpha ) =C+\nu (\alpha ) - \alpha$ and $\s \tau (C) = \s \tau (\rho^{-1}) (\rho C+\s (\nu (\alpha ) - \alpha )).$ $\Box$

{\bf  Proof of Theorem \ref{B6Apr}}. By Theorem \ref{2Apr}, $E= \CD [x,y;\s , \tau , h]$ is a GWA and the element $h$ is a left regular element of $E$. Hence, the element $C= h+\alpha$ is a left regular element of $E$.

$(1\Leftrightarrow 2)$ The equivalence follows at once from Theorem \ref{A6Apr}. $\Box$.

{\bf The canonical central  element $C$ of a diskew polynomial ring (under certain conditions)}.
The next corollary is a criterion for an element $C+\alpha$ (where $\alpha \in D$) to be a central element in $E$. It follows straightaway from Theorem \ref{B6Apr}. This is a generalization of a similar result for the rings $D\langle \s, b , \rho \rangle$,  and {\cite[Lemma 1.5]{Bav-GlGWA-1996}}.

\begin{corollary}\label{aB6Apr}%\marginpar{aB6Apr}
Let $E= D[x,y; \s , \tau , b, \rho ]$ be a diskew polynomial ring such that $\rho$ is a unit, $\CD = D[h; \nu = \tau \s ]$ and $C= h+\alpha$ where $h=yx$ and $\alpha \in D$.  The following statements are equivalent.
\begin{enumerate}
\item The element $C$ is a central element of $E$.

\item $\rho =1$, $\nu =1$, $\alpha - \s (\alpha ) = b$, and the element $\alpha$  belongs to the centre of $D$.

If one of the equivalent conditions holds then
\begin{enumerate}
\item  $C=xy+\s (\alpha )$.
\item $E\simeq D[C] [x,y; \s , \tau , a:= C-\alpha ]$ is a GWA where $\s (C) = C$ and $\tau (C)=C$.
\item The element $C$ is a  regular element of $E$.
\end{enumerate}
\end{enumerate}
\end{corollary}

{\em Every simple ring is, in fact, an algebra}  either over the field of rational numbers $\Q$ or over the finite prime field $\mF_p$ that contains $p$ elements ($p$ is a prime number).

{\bf Simplicity criterion for DPRs in characteristic zero.} {\bf  Proof of Theorem \ref{QDi31Oct}}. $(1\Rightarrow 2)$  This implication follows from Theorem \ref{Di31Oct} and Theorem \ref{B6Apr}.

$(2\Leftarrow 1)$ We have to show that the case (c) of Theorem \ref{Di31Oct} holds. Suppose that this is not the case and $p$ is an element that satisfies the conditions (i)--(iii) in the condition (c) of Theorem \ref{Di31Oct}, we seek a contradiction. Since $D$ is a $\Q$-algebra, the equality (\ref{hpdn6}) can be written as $\rho n^{-1}d_{n-1} - \s (n^{-1}d_{n-1}) = b$. Recall that $d_{n-1}d= \nu (d) d_{n-1}$ for all $d\in D$ and $\nu (d_{n-1})=d_{n-1}$. Now, by Theorem \ref{B6Apr}, the element $C$ is a left normal element of $E$ which is not a unit, and so the ring $E$ is not a simple ring, a contradiction. $\Box $

{\bf Simplicity criterion for DPRs in prime characteristic.} Let $p$ be a prime number and $\mF_p = \Z /p\Z$.  Each natural number $n$ can be written uniquely as a finite sum $n=\sum n_ip^i$ where $0\leq n_i<p$, the $p$-{\em adic form} of $n$. For a natural number $n\neq 0$, let $v_p(n) = \min \{ i\, | \, n_i\neq 0\}$. Then $n=n'p^{v_p(n)}$ for some natural number $n'$ such that $p\nmid n'$. Let $m = \sum m_ip^i$ be the $p$-adic form of a natural number $m$. If $n\geq m$ then
%\marginpar{binmFp}
\begin{equation}\label{binmFp}
{n\choose m} = \prod_i {n_i\choose m_i}\;\; {\rm in}\;\; \mF_p,
\end{equation}
see {\cite[Eqn.  (7)]{Bav-DimMultHolModCharp-2009}} (where ${n_i\choose m_i}=0$ if $n_i<m_i$). For all $i,j$ such that $0\leq j\leq i<p$,  ${n_i\choose m_i}\neq 0$ in $\mF_p$. In particular, (for $n\geq m$),

1. ${n\choose m}=0$ iff there is $i\geq 0$ such that ${n_i\choose m_i}=0$; equivalently, ${n\choose m}\neq 0$ iff $m_i\leq n_i$ for all $i\geq 0$.

2. ${np^i\choose mp^i}={n\choose m}$ for all $i\geq 0$.

{\it Definition.} For a nonzero natural number $n$ such that $p\mid n$, a unique number $\tn$ such that  $ 0\leq \tn < n, \;  {n\choose \tn}\neq 0$  and $ {n\choose i}=0$ (in $\mF_p$) for all $\tn <i<n$ is called the $p$-{\em adic neighbour} of $n$. If $n=n_ip^i+n_{i-1}p^{i-1}+\cdots + n_jp^j$ is the $p$-adic form of $n$ then $\tn =n_ip^i+n_{i-1}p^{i-1}+\cdots + (n_j-1)p^j$ is the $p$-adic form of $\tn$. In particular, $ p\mid\tn$.

Let $\CP = \{ p^i\, | \, i\geq 0\}$. For all natural numbers $i\geq 1$, $\widetilde{p^i}=0$, that is ${p^i\choose j}=0$ for all numbers $j$ such that $0<j<p^i$.

{\bf  Proof of Theorem \ref{PDi31Oct}.} $(2\Leftrightarrow 3)$ This follows from (\ref{hpdn7}) and the equalities ${p^i\choose j}=0$ in $\mF_p$ for all numbers $0<j<p^i$. Let ($2'$) be statement 2 of Theorem \ref{Di31Oct}.

$(1\Leftrightarrow 2')$ Theorem \ref{Di31Oct}.

$(1\Rightarrow 2)$ Clearly, $(2'\Rightarrow 2)$. Then $(1\Rightarrow 2)$ since $(1\Leftrightarrow 2')$.

$(3\Rightarrow 2')$ Suppose that the implication $(3\Rightarrow 2')$ is wrong, we seek a contradiction. Then, by Theorem \ref{Di31Oct}, there is an element $p'=h^n+\sum_{i=0}^{n-1}d_ih^i $ that satisfies the conditions (i)-(iii) of statement $2'$.

(S1) $p\mid n$: By (\ref{hpdn6}), $\rho d_{n-1}-\s (d_{n-1})=nb$, and so $n=0$ in $\mF_p$ (by the condition (c)), i.e.  $p\mid n$.

Let $d_n=1$ and $\Supp (p') :=\{ i\, | \, d_i\neq 0\}$, the {\em support} of the element $p'$. Notice that $n\in \Supp (p')$ and $\Supp (p')\neq \{ n\}$ since otherwise, $p'=h^n$, and so $\beta = \rho^{-1}b=0$ (by the condition (ii) of Theorem \ref{Di31Oct}.(2)) but $b\neq 0$ (the condition (d)), a contradiction.

Recall that $\CP = \{ p^i\, | \, i\geq 0\}$. Then $\Supp (p') = S\coprod T$ where $S=\Supp (p')\backslash \CP$ and $T=\Supp (p')\cap \CP$.

(S2) $S\neq \emptyset$: Since otherwise $p'=h^{p^i}+\sum_{j=0}^{i-1}\alpha_jh^{p^j}+\alpha$, a contradiction.

Let $S= \{ n_1>n_2>\cdots >n_t\}$. Notice that ${t\choose i}=0$ for all $t\in T$ and $0<i<t$. So, {\em for all $j$ such that $0<j<n$ we can ignore the terms in (\ref{hpdn7}) corresponding to the elements} $t\in T$ (since the corresponding binomials are equal to zero).

(S3) {\em If $n>n_1$ then $\s (d_{n_1})=\rho^\nu_{n-n_1}d_{n_1} $ and $d_{n_1}$ is a unit in $D$}:  By (\ref{hpdn7}), for $j=n_1$, $\rho^\nu_{n-n_1}d_{n_1}-\s (d_{n_1})=0$ (all the binomials are equal to zero by the choice of $d_{n_1}$). Since the ring $D$ is a $\s$-simple ring and $d_{n_1}\neq 0$, then the element $d_{n_1}$ must be a unit.

(S4) $n_1=n_1'p^m$ {\em for some $m\geq 1$ and a positive integer $n_1'$ not divisible by $p$}: If $n=n_1$ then the result follows from the statement (S1). Suppose that $n>n_1$. Then for $j=n_1-1$ the equality (\ref{hpdn7}) takes the form $\rho^\nu_{n-n_1+1}d_{n_1-1}-\s (d_{n_1-1}) = n_1\s (d_{n_1}) b$. Then $\s (d_{n_1})^{-1} \rho^\nu_{n-n_1+1}d_{n_1-1}= (\rho^\nu_{n-n_1}d_{n_1})^{-1} \rho^\nu_{n-n_1+1}d_{n_1-1}= d_{n_1}^{-1} (\rho^\nu_{n-n_1})^{-1}\rho^\nu_{n-n_1+1}d_{n_1-1} = d_{n_1}^{-1} \nu^{n-n_1} (\rho ) d_{n_1-1}= \nu^{-(n-n_1)}\nu^{n-n_1}(\rho ) $ $d_{n-1}^{-1}d_{n_1-1}$ $= \rho d_{n_1}^{-1}d_{n_1-1}$. Suppose that $p\nmid n_1$, i.e. $n_1\neq 0$ in $\mF_p$. Then multiplying the equality above by the element $n_1^{-1}\s (d_{n_1})^{-1}$ on the left we have the equality $\rho \alpha - \s (\alpha ) = b$ where $\alpha = n_1^{-1}\s (d_{n_1}^{-1}d_{n_1-1})\in D$, $\nu (\alpha ) = \alpha$ and $\alpha d = \nu^{-(n-n_1)+n-n_1+1}(d) \alpha = \nu (d) \alpha$, a contradiction. Therefore, $p\mid n_1$ and the statement (S4) follows.

(S5) $\tn_1\geq p$ {\em and} $p\mid \tn_1$: Since $p\mid n_1$, we must have $p\mid \tn_1$. Since $n_1\not\in \CP$, we must have $\tn_1\geq p$. The following claim is the essence of the proof.  \\

{\sc Claim (`Meeting the $p$-neighbour method')):} {\em For the element $p'$, either $n_1>n_2=\tn_1$ or, otherwise, there is an index $i\geq 2$ such that $n_1>n_2>\cdots >n_i>n_{i+1}=\tn_i\geq \tn_{i-1}\geq \cdots \geq \tn_1$ (where $\tn_i$ is the $p$-neighbour of $n_i$).}\\

Using the Claim we can get a required contradiction. Suppose that $n_1>n_2=\tn_1$. Then, by (S4), $n_1-\tn_1=p^m$ (recall that $m\geq 1$). Suppose that $n=n_1$. By (\ref{hpdn7}), for $j=\tn_1$, $\rho^\nu_{p^m}d_{\tn_1}-\s ( d_{\tn_1})= {n_1\choose \tn_1} b^{p^m}$. Notice that ${n_1\choose \tn_1}\neq 0$ in $\mF_p$. Let $\alpha = {n_1\choose \tn_1}^{-1}d_{\tn}$. Then $\rho^\nu_{p^m}\alpha -\s (\alpha ) = b^{p^m}$, $\nu (\alpha ) = \alpha$ and $\alpha d = \nu^{p^m}(d)\alpha$ for all $d\in D$,  a contradiction.

Suppose that $n>n_1$.  By (\ref{hpdn7}), for $j=n_1$, $\rho^\nu_{n-n_1}d_{n_1}-\s (d_{n_1})=0$, that is $\s (d_{n_1}) = \rho^\nu_{n-n_1}d_{n_1}$. The ring $D$ is $\s$-simple and $d_{n_1}\neq 0$, so the element $d_{n_1}$ is a unit of $D$. By (\ref{hpdn7}), for $j=\tn_1 (=n_2)$, $\rho^\nu_{n-\tn_1}d_{\tn_1}-\s (d_{\tn_1})= {n_1\choose \tn_1} \s (d_{n_1})b^{p^m}$. By multiplying the equality by the element ${n_1\choose \tn_1}^{-1} \s (d_{n_1})^{-1}$ on the left we have the equality $\rho^\nu_{p^m}\alpha - \s (\alpha ) = b^{p^m}$ where $\alpha = {n_1\choose \tn_1}^{-1} \s (d_{n_1}^{-1}d_{\tn_1})$. We have used the following fact:  $\s (d_{n_1})^{-1} \rho^\nu_{n-\tn_1}d_{\tn_1}= (\rho^\nu_{n-n_1}d_{n_1})^{-1} \rho^\nu_{n-\tn_1}d_{\tn_1}= d_{n_1}^{-1} (\rho^\nu_{n-n_1})^{-1}\rho^\nu_{n-\tn_1}d_{\tn_1} = d_{n_1}^{-1} \nu^{n-n_1} (\rho^\nu_{n_1-\tn_1} ) d_{\tn_1}= \nu^{-(n-n_1)}\nu^{n-n_1}(\rho^\nu_{p^m} ) d_{n_1}^{-1}d_{\tn_1} = \rho^\nu_{p^m} d_{n_1}^{-1}d_{\tn_1}$ since $n_1-\tn_1=p^m$. Notice that $\nu (\alpha ) = \alpha$ (since $\nu (d_{n_1})=d_{n_1}$ and $\nu (d_{\tn_1})=d_{\tn_1}$ as $\tn_1 = n_2$) and $\alpha d = \nu^{p^m}(d) \alpha$ for all elements $d\in D$ (since $d_{n_1}^{-1} d_{\tn_1}d= \nu^{-(n-n_1) +n-\tn_1} (d)d_{n_1}^{-1} d_{\tn_1}= \nu^{p^m} (d) d_{n_1}^{-1} d_{\tn_1}$). So, we obtain a contradiction.

 Suppose that the second case of the Claim holds. For $j=n_i$, the equality (\ref{hpdn7}) takes the form
 $\rho^\nu_{n-n_i}d_{n_i}-\s (d_{n_i}) ={n\choose n_i}b^{n-n_i}+ \sum_{j=1}^{i-1}{n_j\choose n_i} \s (d_{n_j})b^{n_j-n_i}=0$, i.e. $\s (d_{n_1}) = \rho^\nu_{n-n_i}d_{n_i}$, since all the binomials in the equality above are equal to zero as $n_1>n_2>\cdots >n_i>n_{i+1}=\tn_i\geq \tn_{i-1}\geq \cdots \geq \tn_1$. The algebra $D$ is $\s$-simple and $d_{n_i}\neq 0$, hence the element $d_{n_1}$ is a unit in $D$. For $j= n_{i+1}=\tn_i$, the equality (\ref{hpdn7}) can be written as follows
  $\rho^\nu_{n-\tn_i}d_{\tn_i}-\s (d_{\tn_i}) = \sum_{j=1}^{i-1}{n_j\choose n_{i+1}} \s (d_{n_j})b^{n_j-n_{i+1}}+{n_i\choose \tn_i} \s (d_{n_i})b^{n_i-\tn_i}={n_i\choose \tn_i} \s (d_{n_i})b^{n_i-\tn_i}$ since all the binomials but the last one are equal to zero as
 $n_1>n_2>\cdots >n_i>n_{i+1}=\tn_i\geq \tn_{i-1}\geq \cdots \geq \tn_1$. By multiplying the equality above by the element ${n_i\choose \tn_i}^{-1} \s (d_{n_i})^{-1}$ on the left we have the equality
 $\rho^\nu_{p^{v_i}}\alpha - \s (\alpha ) = b^{p^{v_i}}$ where $v_i = v_p(n_i)\geq 1$ (since $p\mid n_i$) and $\alpha = {n_i\choose \tn_i}^{-1} \s (d_{n_i}^{-1}d_{\tn_i})$.  We have used the following fact: $\s (d_{n_i})^{-1} \rho^\nu_{n-\tn_i}d_{\tn_i}= (\rho^\nu_{n-n_i}d_{n_i})^{-1} \rho^\nu_{n-\tn_i}d_{\tn_i}= d_{n_i}^{-1} (\rho^\nu_{n-n_i})^{-1}\rho^\nu_{n-\tn_i}d_{\tn_i} = d_{n_i}^{-1} \nu^{n-n_i} (\rho^\nu_{n_i-\tn_i} ) d_{\tn_i}= \nu^{-(n-n_i)}\nu^{n-n_i}(\rho^\nu_{p^m} ) d_{n_1}^{-1}d_{\tn_1}
   = \rho^\nu_{p^{v_i}} d_{n_i}^{-1}d_{\tn_i}$ since $n_i-\tn_i=p^{v_i}$. Notice that $\nu (\alpha ) = \alpha$ (since $\nu (d_{n_i})=d_{n_i}$ and $\nu (d_{\tn_i})=d_{\tn_i}$) and $\alpha d = \nu^{p^{v_i}}(d) \alpha$ for all elements $d\in D$ (since $d_{n_i}^{-1} d_{\tn_i}d= \nu^{-(n-n_i) +n-\tn_i} (d)d_{n_i}^{-1} d_{\tn_i}= \nu^{p^{v_i}} (d) d_{n_i}^{-1} d_{\tn_i}$). So, we obtain a contradiction.

 {\sc Proof of the Claim.} (S6) $S\neq \{ n_1\}$: Since otherwise, by the equality  (\ref{hpdn7}) for $j=n_1$ we would have $\rho^\nu_{n-n_1}d_{n_1}-\s (d_{n_1})=0$, and so $d_{n_1}$ is a unit, and for $j=\tn_1$, $\rho^\nu_{n-\tn_1}d_{\tn_1}-\s (d_{\tn_1})={n_1\choose \tn_1} \s ( d_{n_1}) b^{p^{n_1-\tn_1}}\neq 0$. Hence, $d_{\tn_1}\neq 0$, i.e., $\tn_1\in T$, and we get  a contradiction by repeating the same arguments as above.

  So, either, $n_1>n_2=\tn_1$ or not. In the second case we must have $n_2>\tn_1$, since otherwise by repeating the arguments of the proof of (S6) for $j=n_1,\tn_1$, we would have a contradiction. Hence, ${n_1\choose n_2}=0$ in $\mF_p$, and so $\rho^\nu_{n-n_2}d_{n_2}-\s (d_{n_2}) =0$, by (\ref{hpdn7}) for $j=n_2$, i.e. $\s (d_{n_2}) =\rho^\nu_{n-n_2}d_2$. Therefore, $d_{n_2}$ is a unit of $D$ since $d_{n_2}\neq 0$ and $D$ is a $\s$-simple ring.

 (S7) $p\mid n_2$: By  (\ref{hpdn7}), for $j = n_2-1$, $\rho^\nu_{n-n_2+1}d_{n_2-1}-\s (d_{n_2-1}) = {n_2\choose n_2-1} \s (d_{n_2})b=n_2\s (d_{n_2}) b$. If $p\not| \, n_2$ then repeating the argument of the proof of (S4) we get a contradiction, So, we must have $n_2=0$ in $\mF_p$, i.e. $p\mid n_2$.

 Let $n_1=a_lp^l+a_{l-1}p^{l-1}+\cdots + a_mp^m$ be the $p$-adic form of $n_1$ where $a_l\neq 0$ and $a_m\neq 0$. Then $\tn_1= a_lp^l+a_{l-1}p^{l-1}+\cdots + (a_m-1)p^m$. The inequalities $n_1>n_2>\tn_1$ imply that $p^m = n_1-\tn_1>n_2-\tn_1>0$. Hence, $n_2= a_lp^l+\cdots + (a_m-1)p^m+a_{m-1}p^{m-1}+\cdots +a_{m_2}p^{m_2}$ is the $p$-adic form of $n_2$ where $a_{m_2}\neq 0$ and $m_2\geq 1$ (since $p\mid n_2$). Then $\tn_2= a_lp^l+\cdots + (a_{m_2}-1)p^{m_2}$ and $\tn_2\geq \tn_1$.

 Suppose that we have already constructed elements $n_1>n_2>\cdots >n_i>\tn_i\geq \tn_{i-1}\geq \cdots \geq \tn_1$ such that $p\mid n_j$ for $j=1, \ldots , i$.

 (S8) $S\neq \{ n_1, \ldots , n_i\}$: Otherwise, by (\ref{hpdn7}), for $j=n_i, \tn_i$, we get a contradiction by the same arguments as above. Then either $n_{i+1} = \tn_i$ and we are done or, otherwise, $n_i>n_{i+1}>\tn_i$. Then repeating the same argument for $n_{i+1}$ as in the case $n_2$ we have that $p\mid n_{i+1}$ and $n_i>n_{i+1}>\tn_{i+1}\geq \tn_i$. This process must stop, say on $j$'th step, that is $n_1>n_2>\cdots >n_j>n_{j+1}=\tn_j>\tn_{j-1}\geq \cdots \geq \tn_1$. The proof of the Claim is complete.    $\Box $

Simplicity criteria for ambiskew polynomial rings are given in \cite{Jordan-Wells-2013}.

{\bf Iterated di-skew polynomial rings.}

  {\em Definition.} A ring of the type $E=D[x_1,y_1;\s_1, \tau_1, b_1, \rho_1]\cdots [x_n,y_n;\s_n, \tau_n, b_n, \rho_n]$ is called an {\bf iterated di-skew polynomial ring} of rank $n$.
 \begin{corollary}\label{a3Jul}%\marginpar{a3Jul}
Let $E=D[x_1,y_1;\s_1, \tau_1, b_1, \rho_1]\cdots [x_n,y_n;\s_n, \tau_n, b_n, \rho_n]$ be an  iterated di-skew polynomial ring of rank $n$. Then $E=\oplus_{\alpha , \beta\in \N^n}D(yx)^{\alpha , \beta}$ is a free left $D$-module where $(yx)^{\alpha , \beta}:=y_1^{\alpha_1}x_1^{\beta_1}y_2^{\alpha_2}x_2^{\beta_2}\cdots y_n^{\alpha_n}x_n^{\beta_n}$, $\alpha = (\alpha_1, \ldots , \alpha_n)$ and $\beta = (\beta_1, \ldots , \beta_n)$.
\end{corollary}

{\it Proof}. The corollary follows from Theorem \ref{5Apr}. $\Box $

\small{

Department of Pure Mathematics

University of Sheffield

Hicks Building

Sheffield S3 7RH

UK

email: v.bavula@sheffield.ac.uk
}

\end{document}